\documentclass[11pt]{article}
\usepackage[utf8]{inputenc} 
\usepackage[T1]{fontenc}	
\usepackage{lmodern}		
\usepackage[british]{babel}	
\usepackage[a4paper, total={160mm, 240mm}]{geometry}
\usepackage{amsmath}
\usepackage{amssymb}
\usepackage{amsthm}
\usepackage{mathrsfs}
\usepackage{mathtools}
\usepackage{minted}
\usepackage{bm}
\usepackage{enumerate}
\usepackage[shortlabels]{enumitem}
\usepackage{tikz-cd}
\usepackage{comment}
\usepackage[thinc]{esdiff}
\usepackage{hyperref}
\usepackage[style=numeric,maxbibnames=99]{biblatex}
\usepackage{array}

\usepackage{booktabs}
\usepackage{tabularx}
\usepackage{makecell}
\newcommand{\tabitem}{~~\llap{\textbullet}~~}

\definecolor{pastel1}{RGB}{245, 235, 235}
\definecolor{pastel2}{RGB}{238, 241, 255}
\definecolor{pastel3}{RGB}{255, 249, 249}
\definecolor{pastel4}{RGB}{247, 232, 246}
\definecolor{pastel5}{RGB}{253, 250, 246}
\definecolor{pastel6}{RGB}{228, 239, 231}
\definecolor{pastel7}{RGB}{243, 241, 245}

\hypersetup{
	colorlinks=true,
	linkcolor=blue,
	citecolor=red,
}

\newtheorem{theorem}{Theorem}[section]
\newtheorem{prop}[theorem]{Proposition}
\newtheorem{lemma}[theorem]{Lemma}

\theoremstyle{definition}
\newtheorem{definition}[theorem]{Definition}
\newtheorem{example}[theorem]{Example}

\newtheorem{classification}[theorem]{Classification}
\newtheorem{corrections}[theorem]{Correction}
\newtheorem{obs}[theorem]{Remark}

\theoremstyle{remark}
\newtheorem{question}[theorem]{Question}

\DeclareMathOperator{\Ker}{Ker}

\DeclareMathOperator{\Ind}{Ind}
\DeclareMathOperator{\Res}{Res}

\newcommand{\MSh}{M_{\mathrm{Sh}}}

\newcommand{\C}{\mathbb{C}}

\newcommand{\PGL}[2]{\mathrm{PGL}_{#1}(#2)}
\newcommand{\PSL}[2]{\mathrm{PSL}_{#1}(#2)}
\newcommand{\PSU}[2]{\mathrm{PSU}_{#1}(#2)}
\newcommand{\SU}[2]{\mathrm{SU}_{#1}(#2)}
\newcommand{\Sp}[2]{\mathrm{Sp}_{#1}(#2)}
\newcommand{\PSp}[2]{\mathrm{PSp}_{#1}(#2)}

\newcommand{\GL}[2]{\mathrm{GL}_{#1}({#2})}
\newcommand{\SL}[2]{\mathrm{SL}_{#1}({#2})}

\newcommand{\Kr}{\mathrm{Kr}}
\newcommand{\Pc}{\mathbb{P}}

\newrobustcmd*{\citeallauthors}{%
	\AtNextCite{\AtEachCitekey{\defcounter{maxnames}{999}}}%
	\citeauthor}

\title{On the finite irreducible subgroups of $\GL{n}{\C}$ and $\PGL{n}{\C}$}
\author{Gerard Gonzalo Calbetó}

\begin{document}
	\maketitle
	
	\begin{abstract}
		The classification of the finite subgroups of $\GL{n}{\C}$ and $\PGL{n}{\C}$ is a classical problem in the field of finite group theory, dating back to the late 19th century with authors like Klein, Jordan, Blichfeldt, etc. Throughout its long history, many results concerning the classification have been scattered in the mathematical literature. In this survey, we explore many of the most relevant results related to the classification, its structure, and the lists of groups. We mostly focus on the irreducible groups. In particular, classification statements are provided for primitive and imprimitive groups over prime dimension, and quasi-primitive groups of small composite dimension. We also provide tables with a detailed list of all (quasi)-primitive finite groups of dimension $n<8$. 
		
		Finally, we provide a computer program implementing many of the known results specific to the finite quasisimple irreducible projective groups, to improve and preserve the accessibility to these results and further their classification. See \S \ref{subsec:computational} and \cite{FISGO}.
		
		We aim this survey to the non-specialist, so the provided classification results may easily accessible to mathematicians working outside of group theory.
	\end{abstract}
	
	\tableofcontents
	
	\newpage
	\section{Introduction}
	
	The classification of finite simple groups is, to this date, the most renowned classification theorem in group theory, and possibly one of the greatest undertakings in mathematics. Another classical classification problem, born around the same time as the former, is the classification of the finite groups of the linear complex group $\GL{n}{\C}$ for a fixed dimension $n$. This problem also includes the classification for the unimodular group $\SL{n}{\C}$ and their corresponding projective analogs, $\PGL{n}{\C}$. We remark that classifications of $\GL{n}{\C}$ and $\SL{n}{\C}$ both induce a classification of $\PGL{n}{\C} \cong \PSL{n}{\C}$ via projection, and much of the literature involves either $\GL{n}{\C}$ or $\SL{n}{\C}$.
	
	This classification endeavor dates back to authors like Jordan \cite{Jordan1877}, Klein \cite{Klein1884}, E.H.Moore, Dickson, \citeauthor{Blichfeldt1917} \cite{Blichfeldt1917} and many others, with a vast literature extending from the early 1870's to the current day. In particular, many surveys on the topic can be found such as \cite{Feit1971} by Feit, \cite{Zalesskii1989} by \citeauthor{Zalesskii1989} and \cite{Tiep2000} by \citeauthor{Tiep2000}. As for books specifically devoted to the topic, the main sources are \citeauthor{Blichfeldt1917}'s book \cite{Blichfeldt1917} and \citeauthor{Dixon1971}'s \cite{Dixon1971}.
	
	In this survey, we aim to make available an updated collection of the most relevant information concerning these classifications. We specially gear it towards the non-specialists, so it may be of use to the general mathematician, whichever their area of work. In this regard, it is worth mentioning that most results described over $\C$ can be generalized to algebraically closed fields of zero characteristic and, under certain restrictions, to positive characteristic. Although it will not be reviewed here, the interested reader may refer to \cite[\S 5]{Dixon1971}.
	
	The description of the finite groups of $\PGL{n}{\C}$ starts with Klein \cite{Klein1884}, who described the case $n=2$. The case $n=3$ appears listed in great generality by \citeauthor{Huggins2006} in \cite{Huggins2006}. Recently, \citeauthor{Flannery2021} in \cite{Flannery2021} have provided a complete description of all irreducible groups of $\GL{2}{\C}$ and $\GL{3}{\C}$ alongside a computational implementation \cite{Monomial2021} based on Magma \cite{Magma}.  This allows access to all irreducible representations of finite groups in  $\GL{2}{\C}$ and $\GL{3}{\C}$. 
	
	For $n > 3$, the classifications are usually split between primitive and imprimitive groups, which we introduce in \S \ref{sec:outline}. Case $n=4$ was explored by \citeauthor{Blichfeldt1917} and culminated in his book \cite{Blichfeldt1917} with a full description of the primitive groups of $\SL{4}{\C}$. The explicit description of the imprimitive groups for $\GL{4}{\C}$ is more recent and due to \citeauthor{Flannery1999} in \cite{Flannery1999} and \citeauthor{Hofling2001} in \cite{Hofling2001}. The primitive groups of $\SL{n}{\C}$ were described for $n=5$ by \citeauthor{Brauer1967} in \cite{Brauer1967}, for $n=6$ by \citeauthor{Lin71} in \cite{Lin71}, for $n=7$ by \citeauthor{Wales1970} in \cite{Wales1970}, for $n=8$ by \citeauthor{Feit1976}, \citeauthor{Huffman1976} in \cite{Feit1976} and \cite{Huffman1976} and finally for $n=11$ by \citeauthor{Rob81} in \cite{Rob81}. The classifications of $n=9,10$ currently contain gaps, we go into more detail in section \S\ref{sec:composite}.
	
	In a more general context, relatively recent developments, and the classification of finite simple groups, have enabled an almost complete description of the irreducible groups of $\PGL{p}{\C}$, $\SL{p}{\C}$ and $\GL{p}{\C}$ when the dimension $p$ is a prime number. The descriptions of the imprimitive case over $\GL{p}{\C}$ are due to \citeauthor{Flannery2021} in \cite{Flannery2021}, who also provide a computational implementation of their results in \cite{Monomial2021}, and due to \citeauthor{DixonZalesskii2004Imprimitive} in \cite{DixonZalesskii2004Imprimitive} when restricting to non-solvable groups and over $\SL{p}{\C}$. The description of the primitive case over $\SL{p}{\C}$ is mostly due to \citeauthor{DixonZalesskii1998Primitive} in \cite{DixonZalesskii1998Primitive} and \cite{DixonZalesskii2008Corrigendum}, and due to \citeallauthors{KangZhang2009} in \cite{KangZhang2009}.
	
	Before proceeding to a description of current applications of these results and the structure of this survey, for the reader's convenience we provide a collection of the most relevant the classification statements included in this survey and where to find them.
	
	\begin{itemize}
		\item Imprimitive groups over prime dimension: Section~\ref{subsec:structure_imprimitive}, Classification~\ref{class:monomial}.
		\item Primitive groups over small prime dimension $2\leq p \leq 11$: Section~\ref{subsubsec:classic_approach}, Classification~\ref{class:primitive_small_p}.
		\item Primitive groups over prime dimension: Section~\ref{subsubsec:post_approach}, Classification~\ref{class:primitive_any_p}.
		\item Primitive groups over small composite dimension $n < 10$: Section~\ref{sec:composite}, Classification~\ref{class:primitive_small_composite}.
		\item Minimal degree of an irreducible projective simple group of Lie type in zero characteristic: Section~\ref{subsec:quasisimple}, Classification~\ref{class:small_pirreps}.
		\item Irreducible (quasi-primitive) projective simple groups of dimension up to $250$: Section~\ref{subsec:quasisimple}, Classification~\ref{class:HissMalle}.
	\end{itemize}
	
	Currently, we are developing a program based on Python that will allow users to compute a finite list of simple groups candidates to being quasi-primitive projective groups of dimension $n$, alongside information on the existence of representations for each candidate. Our main objective is to provide a computational interface to many of the long classification results available as tables in the literature, many of them present in this survey. The details are discussed in Section~\ref{subsec:computational}, furthermore, development progress and the program itself can be consulted in the author's GitHub \cite{FISGO}. \\
	
	Recently there has been a renewed interest in a complete list of groups of the projective linear spaces in fields such as algebraic geometry or arithmetic geometry. 
	
	In algebraic geometry, given a smooth projective variety $\overline{V}$ in $\Pc^n_{\overline{k}}$ over an algebraically closed field $\overline{k}$, the group of birational transformations of $\overline{V}$ onto itself, denoted by $\mathrm{Bir}(\overline{V})$, contains the group of its automorphisms $\mathrm{Aut}(\overline{V})$. Moreover, we denote by $\mathrm{Lin}(\overline{V}) \leq \mathrm{Aut}(\overline{V})$ its subgroup of automorphisms induced by projective linear transformations $\mathrm{Aut}(\Pc^n_{\overline{k}}) = \PGL{n+1}{\overline{k}}$. A classical problem, which goes back to Cremona \cite{Cremona1863} in 1863, is the study of $\mathrm{Bir}(\Pc^2_{\overline{k}})$ and its relation to $\mathrm{Lin}(\Pc^2_{\overline{k}}) = \PGL{3}{\overline{k}}$. Apart from algebraic geometry, this problem is also a subject of interest in other fields such as dynamical systems. For instance, see results of \citeauthor{Blanc2021} in \cite{Blanc2021} or Dolgachev in \cite{Dolgachev2010} among many others, which we apologize for omitting.
	
	Another field within algebraic geometry where the finite subgroups of $\SL{n}{\C}$ are of special relevance is in the study of quotient singularities. This is explained in great detail in the introduction of \cite{KangZhang2009}, so we refer the interested reader to it.
	
	In arithmetic geometry, smooth projective varieties with non-trivial automorphisms are usually of key interest. Examples of such curves are Fermat's or Klein's degree $d\geq 4$ curves over $\Pc^2_k$. Furthermore, the study of the moduli of smooth plane curves with certain automorphism groups are also currently of interest, see for instance \cite{Badr2016}, \cite{Badr2022} or \cite{Lercier2014}. In particular, the automorphism group of a smooth plane curve is a finite subgroup of $\PGL{3}{\overline{k}}$ which fixes the equation of the curve, so a list of the finite subgroups of $\PGL{3}{\overline{k}}$ is of key interest.
	
	In a more general setting, let $\overline{V}$ be a $(n-1)$-dimensional smooth projective variety identified with a hypersurface model $H_{\overline{V},d,n}$ represented by a single homogeneous polynomial equation without singularities. Let $F(X_0,\ldots,X_n)=0$ be such an equation for some degree $d$ over $\overline{k}$ (assume once and for all that $d\geq 4$), where $X_0,\ldots,X_n$ denote the homogeneous coordinate system for the projective space. In 1964, Matsumura and Monsky \cite{Matsumura1964} showed that, for $n\geq 3$, $\mathrm{Lin}(\overline{V})$ is a finite group and $\mathrm{Aut}(\overline{V})=\mathrm{Lin}(\overline{V})\leq \mathrm{PGL}_n(\overline{k})$ except possibly when $(n,d)=(3,4)$. In this direction, results of \citeauthor{Blichfeldt1917} about $\PGL{4}{\C}$ have recently been used by \citeauthor{Cheltsov2019} in \cite{Cheltsov2019} and by \citeauthor{Avila2024} in \cite{Avila2024}. \\
	
	We now detail section by section the contents presented in this survey. Section~\ref{sec:structure_section} is mostly devoted to introducing the structure of the classification. Section~\ref{subsec:notation} includes notation conventions to be used throughout the survey. In Section~\ref{sec:prelims_proj_reps} we describe the interplay between $\GL{n}{\C}$, $\SL{n}{\C}$ and $\PGL{n}{\C}$ giving an introduction to the theory of projective representations. Section~\ref{sec:outline} provides definitions of the different types of groups used to split the classification alongside examples. Section~\ref{subsec:structure_imprimitive} introduces the main structure of imprimitive groups and the classification results regarding imprimitive groups, contained in Classification~\ref{class:monomial}. In section~\ref{subsec:structure_primitive} we introduce basic structural results regarding primitive and quasi-primitive groups alongside many bounds for the order of these groups and the primes dividing them.
	
	We devote Section~\ref{sec:primitive_groups} to describing the classifications of primitive groups in $\SL{n}{\C}$ and $\PGL{n}{\C}$, which we split into two parts. Firstly, Section~\ref{sec:primitive_prime} tackles the classification over prime dimension, presenting both the classical approach and the modern (post-classification of finite simple groups) approach. In this section we include Classification~\ref{class:primitive_small_p} of primitive groups of small prime dimension, alongside Corrections~\ref{correction:small_p_1}, \ref{correction:small_p_2} mending some errors contained in the classical results; the section ends with Classification~\ref{class:primitive_any_p} of the primitive groups over any prime dimension. Both classification statements are relevant, for the latter is less explicit than the former. Finally, Section~\ref{sec:composite} covers the case of composite dimension, presenting an important result due to \citeauthor{Lin69} on their structure, Classification~\ref{class:primitive_small_composite} covering the known classifications of the primitive groups over small composite dimension, and finally present Theorem~\ref{teoTZ96} describing all irreducible representations of quasisimple groups for relatively small dimension.
	
	Finally, Section~\ref{sec:tables} mainly contains practical information in the form of Tables~\ref{tab:pgl_2}, \ref{tab:pgl_3}, \ref{tab:pgl_4}, \ref{tab:pgl_5}, \ref{tab:pgl_6}, \ref{tab:pgl_7} listing all the primitive groups in dimensions $2\leq n \leq 7$. Furthermore, Table~\ref{tab:missing} contains a status of the classification of primitive groups of degree $2\leq n \leq 11$ with comments on their completion or missing parts. Section~\ref{subsec:computational} contains a description of the software \cite{FISGO}, currently in development by the author, to possibly determine the quasi-primitive simple projective groups for any dimension.\\
	
	\vspace*{\fill}

	\textbf{Acknowledgments} The author is very grateful to Francesc Bars, who directed the author's bachelor thesis on this classification problem and encouraged its refinement in the form of this survey. His support and encouragement during the preparation of this survey, alongside his many comments and revisions, have proven to be invaluable. 
	
	We would also like to thank professor G.R. Robinson for his assistance in accessing his PhD thesis, and all the valuable comments and feedback on the earlier versions of this survey. 
	
	\newpage
	\section{Structure of finite linear and projective groups} \label{sec:structure_section}
	
	Throughout this survey, we focus our efforts on understanding the classification of finite subgroups of the complex projective general linear group $\PGL{n}{\C},\, n\geq 2$ with fixed $n$. We refer to the quantity $n$ as either \emph{dimension} or \emph{degree}. We also are required to talk about its linear counterparts, mainly the general linear group $\GL{n}{\C}$ and its subgroup of matrices with unit determinant $\SL{n}{\C}$. Also recall that $\PSL{n}{\C} = \PGL{n}{\C}$ \cite[\S 2.2]{Jones1987}.
	
	\subsection{Notation} \label{subsec:notation}
	
	We start fixing some basic notions and notation. We denote by $G$ a finite abstract group, $Z(G)$ its center, $G' = [G,G]$ its derived or commutator subgroup, $\mathrm{Aut}(G)$ its group of automorphisms, $\mathrm{Inn}(G)$ the group of inner automorphisms (i.e. automorphisms given by conjugation by a fixed element of $G$), and $\mathrm{Out}(G) = \mathrm{Aut}(G)/\mathrm{Inn}(G)$ the outer automorphism group. Given a subgroup $H\leq G$ we denote by $[G:H]$ its index inside $G$.
	
	A finite subgroup $H\leq\GL{n}{\C}$ will be denoted by $H = \rho(G)$, where $G$ is the abstract group isomorphic to the linear group $H$ and $\rho: G\to \GL{n}{\C}$ its corresponding faithful linear representation. Similarly, we denote projective subgroups $H\leq \PGL{n}{\C}$ as $H = \bar{\rho}(G)$ where $\bar{\rho}:G \to \PGL{n}{\C}$ is the corresponding faithful projective representation. We denote the Kronecker product of matrices or representations by $\otimes_{\Kr}$.
	
	We also introduce notation to refer to certain types of finite groups.
	\begin{itemize}
		\item Cyclic groups of order $n$: $C_n$.
		\item Alternation group on $n$ points: $A_n$.
		\item Symmetric group on $n$ points: $S_n$.
		\item Projective special group of dimension $n$ over a finite field $\mathbb{F}_q$: $\PSL{}{n,q}$
	\end{itemize}
	
	Finally, for some groups $A,B,C$, when we write $0\to A\to B \to C\to 0$ we always consider it a short exact sequence unless stated otherwise.
	
	\subsection{Preliminaries on projective representations} \label{sec:prelims_proj_reps}
	This subsection serves as a primer on projective representation theory. A reader familiar with the theory may skip to Section \ref{sec:outline}. We mostly introduce necessary concepts to understand how the classifications of $\GL{n}{\C}$, $\SL{n}{\C}$ and $\PGL{n}{\C}$ interact. Most of the coming exposition is covered in \cite[\S 3]{Karpilovsky1985} and \cite[\S 1]{Hoffman1992}.\\
	
	To understand the interplay between group representations in $\GL{n}{\C}$ and $\PGL{n}{\C}$, it is convenient to introduce the following notion for projective representations.
	\begin{definition} \label{def:projective_rep}
		A map $\tau: G \to \GL{n}{\C}$ is a projective representation of $G$ over $\C$ if there exists a map $\alpha : G\times G \to \C^\times$ such that 
		\begin{enumerate}[(i)]
			\item $\tau(x)\tau(y) = \alpha(x,y)\tau(xy)$ for all $x,y\in G$.
			\item $\tau(1) = 1$.
		\end{enumerate} 
		In particular, if $\pi: \GL{n}{\C} \to \PGL{n}{\C}$ is the natural projection, then $\bar{\tau} = \pi\circ\tau : G \to \PGL{n}{\C}$ is a group homomorphism. We will commonly abuse the notation and also refer to $\bar{\tau}$ as a projective representation.
	\end{definition}
	\begin{obs}
		We sometimes denote a projective representation as a pair $(\tau,\alpha)$ if the $\alpha$ map is explicitly needed. It is readily seen that $\alpha$ satisfies the properties of a $2$-cocycle $\alpha\in Z^2(G,\C^\times)$. Furthermore, note that if $\alpha \equiv 1$ then $\tau$ is a linear representation.
	\end{obs}
	From the point of view of projective representations as $\bar{\tau}:G \to \PGL{n}{G}$ group homomorphisms, conjugation is a natural equivalence relation for representations. Meaning $\bar{\tau}$ and $\bar{\tau}'$ are equivalent if and only if there exists $T\in \PGL{n}{\C}$ such that $\bar{\tau}(g) = T^{-1}\tau(g)T$ for all $g\in G$. Thus, we can wonder is the correspondence between pairs $(\tau,\alpha)$ and $\bar{\tau}$ is unique up to conjugation. This is generally false, meaning there might be some $(\rho,\beta) \neq (\tau,\alpha)$ such that $\bar{\rho}$ is conjugate to $\bar{\tau}$. This motivates the following definition.
	\begin{definition}
		Two projective representations $\tau_1, \tau_2$ of degree $n$ are projectively equivalent if there exists a map $\mu:G\to \C^\times$ with $\mu(1)=1$ and $T\in\GL{n}{\C}$ such that $\tau_2(g) = \mu(g)T\tau_1(g)T^{-1}$ for all $g\in G$. 
	\end{definition}
	\begin{example}
		Take $G$ such that there exists a linear representation $\rho: G\to \GL{n}{\C}$ and a non-trivial 1-dimensional representation $\alpha: G \to \C^\times$ such that $\rho' \coloneqq\rho\otimes\alpha$ is a different linear representation from $\rho$. It follows that $\rho$ and $\rho'$ are projectively equivalent. 
	\end{example}
	\begin{obs} 
		Thus we conclude that two projective representations $(\tau,\alpha)$ and $(\rho,\beta)$ are projectively equivalent if and only if their group homomorphisms $\bar{\tau},\,\bar{\rho}$ are conjugate. Furthermore, $2$-cocycles $\alpha,\beta\in Z^2(G,\C^\times)$ produce projectively equivalent representations if and only if $\alpha,\beta$ are cohomologous. So the abelian group $M(G) = H^2(G,\C^\times)$, called the Schur multiplier of $G$, controls the projective representations of $G$. 
	\end{obs} 
	Ideally, we want all projective representations to arise from linear representations. However, the linear representations of the group itself are not enough, take any representation arising from a non-trivial 2-cocycle. This problem can be essentially reduced to studying specific central extensions of $G$. Recall that an extension $0\to A\overset{i}{\to} B \to C \to 0$ is called \emph{central} if $i(A) \subseteq Z(B)$. A central extension is called a \emph{stem} extension if $i(A) \subseteq Z(B) \cap B'$.
	\begin{definition}
		A representation group of $G$, denoted by $G^*$, is a stem extension $0\to \MSh\to G^* \to G\to 0$, i.e. the sequence is exact with $\MSh\cong M(G)$, and $\MSh \subseteq (G^*)'\cap Z(G^*)$.
	\end{definition}
	Depending on the author, representation groups are sometimes referred to as Schur coverings of $G$. This notion is not to be confused with \emph{coverings} of $G$, which refer to extensions of $G$ by a certain group $A$ isomorphic to a subgroup of $M(G)$ (see Definition \ref{def:covering}).
	\begin{obs}
	For any group $G$, there exists a representation group $G^*$ \cite[Theorem 1.2]{Hoffman1992} which need not be unique (see \ref{rmk:unique_rep_group}). Representation groups can be interpreted as liftings of projective representations through the following diagram of exact sequences
		\begin{equation} \label{eq:diagram1}
			\begin{tikzcd}
				0 \arrow[r] & \MSh \arrow[r]\arrow[d, "\beta"] & G^* \arrow[r,"f"]\arrow[d,dashed, "\rho"] & G \arrow[r]\arrow[d,"\bar\rho"] & 0 \\
				0 \arrow[r] & \C^* \arrow[r, "i"] & \GL{n}{\C} \arrow[r, "\pi"] & \PGL{n}{\C} \arrow[r] & 0
			\end{tikzcd}
		\end{equation}
	Note that if we have a linear representation $\Gamma: G^* \to \GL{n}{\C}$ such that $\Gamma(a)$ is a scalar transformation for any $a\in \MSh$, and a set section $\mu:G \to G^*$ of $f$, then $\tau = \Gamma\circ\mu$ defines a projective representation of $G$. 
	\end{obs}
	The following theorem states that all projective representations of a group arise from linear representations of a representation group. See \cite[Theorem 3.3.7]{Karpilovsky1985} and \cite[Theorems 1.3,1,4]{Hoffman1992}.
	\begin{theorem}
		Let $1\to \MSh \to G^* \to G \to 0$ be a central extension with $|G^*|=|G||M(G)|$. The following are equivalent:
		\begin{enumerate}[(i)]
			\item $G^*$ is a representation group of $G$.
			\item Any projective representation $\bar\rho(G): G\to \PGL{n}{\C}$ can be lifted to a linear representation $\rho:G^* \to \GL{n}{\C}$ such that diagram (\ref{eq:diagram1}) commutes.
		\end{enumerate}
	\end{theorem}
	
	\begin{obs} \label{rmk:unique_rep_group}
		Representation groups need not be unique, and Schur \cite[Theorem 3.4.4]{Karpilovsky1985} gave an upper bound on the number of non-isomorphic representation groups. In particular, perfect groups (i.e if $G = G'$) have a unique representation group.
	\end{obs}
	
	This description of projective representations allows the extension of many properties defined for linear representations into the projective setting, so long as they respect projective equivalence. We finish this primer with some comments on coverings and faithful representations.
	
	\begin{obs}
		Representation groups enable the lifting of all projective representations. However, for a specific representation, it might suffice to describe a ``smaller'' central extension by a subgroup of $\MSh$ instead of the full $\MSh$ group.
	\end{obs}
	
	\begin{definition}\label{def:covering}
		Let $G$ be a group, a covering group of $G$ is a central extension $0\to A \to A.G \to G \to 0$ such that $A \leq \MSh \cong M(G)$ and $A \subseteq (A.G)'$. By abuse of notation, we refer to $A.G$ as a covering of $G$.
	\end{definition}
	\begin{obs}
		Note that coverings of $G$ by $A$ need not be unique, however they are usually denoted by $A.G$	in the same way semidirect products are denoted $A\rtimes G$. Furthermore, if there exists a unique covering $A.G$ where $A\cong C_n$ is cyclic, it is usually denoted as $n.G$.
	\end{obs}
	
	\begin{example}
		Consider the alternating group $A_6$. Then $M(A_6) = C_6$, so we may consider the coverings given by $C_1, C_2, C_3$ and $C_6$. For instance, linear representations $\tau: A_6 \to \GL{n}{\C}$ projected onto $\PGL{n}{\C}$, i.e. $\bar{\tau} = \pi\circ\tau: G \to \PGL{n}{\C}$, correspond to projective representations lifted through $C_1$. Projected linear representations of $2.A_6$ (unique covering of $A_6$ by $C_2$) are the projective representations of $A_6$ which can be lifted through $C_2$, these include the ones lifted by $C_1$, and so on.
	\end{example}
	
	\begin{obs} \label{rmk:similar_diagram}
		Let $\bar{\rho}:G \to \PGL{n}{\C}$ be a projective representation that can be lifted through $A \leq \MSh$. Then a covering group $A.G$ produces the same diagram as (\ref{eq:diagram1}) exchanging $\MSh$ by $A$ and $G^*$ by $A.G$. 
	\end{obs}
	
	Recall that a faithful representation is an injective representation. Understanding how faithfulness translates from the projective setting to the linear setting is crucial to understand the interaction between classifications of finite groups in $\GL{n}{\C}$ and $\PGL{n}{\C}$.
	
	\begin{prop} \label{prop:lift_faithful}
		Let $\bar\rho: G\to \PGL{n}{\C}$ be a faithful projective representation, $M(G)$ its Schur multiplier and $G^*$ a representation group. Then there exists a subgroup $A\leq \MSh \cong M(G)$ such that the covering $A.G$ given by $0\to A \to A.G \to G \to 0$ lifts $\bar{\rho}$ into a linear faithful representation of $A.G$.
	\end{prop}
	\begin{proof}
		Let $G^*$ be a representation group of $G$ and consider Diagram~(\ref{eq:diagram1}) with $\rho: G^*\to \GL{n}{\C}$ a lifting of $\bar\rho$ and $\pi:\GL{n}{\C}\to\PGL{n}{\C}$ the canonical projection. Then $ \{1\} = (\pi\circ\rho)(\ker\rho) = (\bar\rho\circ f)(\ker\rho)$, and by injectivity of $\bar{\rho}$ we have $f(\ker\rho) = \{1\}$ so $\ker\rho \leq \MSh$. Thus, we may consider the following commutative diagram 
		\begin{equation}
			\label{eq:diagram2}
			\begin{tikzcd}
				0 \arrow[r] & \MSh/\ker\rho \arrow[r]\arrow[d, "\tilde\beta"] & G^*/\ker\rho \arrow[r,"\tilde f"]\arrow[d,dashed, "\rho'"] & G \arrow[r]\arrow[d,"\bar\rho"] & 0 \\
				0 \arrow[r] & \C^* \arrow[r, "i"] & \GL{n}{\C} \arrow[r, "\pi"] & \PGL{n}{\C} \arrow[r] & 0
			\end{tikzcd}
		\end{equation}
		and obtain a faithful representation $\rho': G^*/\ker\rho \to \GL{n}{\C}$ whose projection is $\bar{\rho}$. Since $\MSh$ is abelian, there exists $A\leq \MSh$ with $A \cong \MSh/\ker\tau$. We need to check that the extension is stem. Since $Z(G^*/\ker\rho) \supseteq Z(G^*)/\ker\rho \supseteq \MSh/\ker\rho$ the extension is central. Furthermore, $(G^*/\ker\rho)' = (G^*)'/\ker\rho$ ($G'$ is defined by words $aba^{-1}b^{-1}$ so it is preserved), so the extension is stem.
	\end{proof}
	\begin{prop} \label{prop:proj_faithful}
		Let $A\leq M(G)$ and consider a covering $A.G$ such that $\rho:A.G \to \GL{n}{\C}$ is faithful. Then the projected representation $\bar\rho = \pi\circ\rho$ is faithful and $A$ is cyclic.
	\end{prop}
	\begin{proof}
		We use Remark~\ref{rmk:similar_diagram} and notation from diagram (\ref{eq:diagram1}). Note that if $\rho$ is injective then $\beta:A\to \C^*$ must also be injective. In particular, $A$ is cyclic. We can restrict the diagram to obtain the commutative diagram
		\begin{equation*}
			\begin{tikzcd}
				0 \arrow[r] & A \arrow[r]\arrow[d,"\beta"] & A.G \arrow[r,"f"]\arrow[d, "\rho"] & G \arrow[r]\arrow[d,"\bar\rho"] & 0 \\
				0 \arrow[r] & \beta(A) \arrow[r] & \rho(A.G) \arrow[r, "\pi|_{\rho(A.G)}"] & \PGL{n}{\C} & 
			\end{tikzcd}
		\end{equation*}
		and by the Snake lemma we conclude $\bar\rho$ is injective.
	\end{proof}
	\begin{obs} \label{rmk:unimodular_lift}
		Propositions~\ref{prop:lift_faithful}, and \ref{prop:proj_faithful} cover how finite groups ascend and descend from $\GL{n}{\C}$ to $\PGL{n}{\C}$ and viceversa. However, many results over the next sections are covered over $\SL{n}{\C}$, so need to know if these lifts can be done in a unimodular way. Lifting finite groups from $\PGL{n}{\C}$ to $\SL{n}{\C}$ is much simpler. Let $\overline{H}\leq \PGL{n}{\C}$ be a finite subgroup, since $\PGL{n}{\C} \cong \PSL{n}{\C}$, we may choose representatives $h \in\SL{n}{\C}$ of each $\overline{h}\in \overline{H}$, with $\pi(h) = \overline{h}$ and $\pi: \SL{n}{\C}\to \PGL{n}{\C}$. Then $H \coloneqq \langle h \mid \overline{h}\in\overline{H} \rangle \leq \SL{n}{\C}$ is a subgroup and $|H|\leq |\overline{H}||Z(\SL{n}{\C})| < \infty$ as $Z(\SL{n}{\C})\cong C_n$ is finite. The resulting group $H$ is such that $\pi(H) = \overline{H}$. Note this same argument does not work over $\GL{n}{\C}$ for $|Z(\GL{n}{\C})| = \infty$ and so the theory of projective representations is needed. 
	\end{obs}

	\subsection{Outline of the classifications} \label{sec:outline}
	Recall a representation $\rho:G \to \GL{n}{\C}$ is said to be irreducible if there is no proper vector space $V \subseteq \C^n$ such that $\rho(G)(V) \subseteq V$, otherwise, it is said to be reducible. Over $\C$, any reducible representation decomposes as the direct sum of irreducible representations (see Maschke's theorem). This is not necessarily the case over general fields, and when it happens the representation is said to be completely reducible. 
	\begin{definition} \label{def:irred}
		An irreducible representation $\rho:G \to \GL{n}{\C}$ is \emph{imprimitive} if there exist proper nonzero vector spaces $V_1, \dots, V_k$ such that $\C^n = V_1 \oplus \cdots \oplus V_k$ and for any $g\in G$ and all $i=1,\dots,k$, $\rho(g)(V_i) = V_j$ for some $j=1,\dots,k, \, j\neq i$ . The set $\{V_1, \dots, V_k\}$ is a \emph{set of imprimitivity} of $\rho(G)$. If such a set does not exist, $\rho$ is \emph{primitive}.
	\end{definition}
	A geometric interpretation of these notions is given in the introduction of \cite{Cheltsov2019} for $n=4$. Furthermore, since the $\rho(g)$ are vector space isomorphisms, the vector spaces $V_1, \dots, V_k$ must all have the same dimension. 
	\begin{definition} \label{def:imprimitive}
		Let $\rho:G \to \GL{n}{\C}$ be a faithful imprimitive representation alongside proper nonzero vector spaces $V_1, \dots, V_k$ whose dimension $d$ is minimal among all possible sets of imprimitivity. We say that $\rho(G)$ is $d$\emph{-nomial}. If $d=1$ we say it is \emph{monomial}.
	\end{definition}
	\begin{obs}
		Definitions~\ref{def:irred} and \ref{def:imprimitive} provide a classical notion of primitivity and imprimitivity. These notions can also be characterized in representation theoretic terms using \emph{induced representations}.
	\end{obs}
	\begin{definition} \label{def:induced_rep}
		Let $G$ be a finite group and $H \leq G$ a subgroup. Consider $k = [G:H]$ the size of the left cosets $G/H$, and $R  = \{r_1, \dots, r_k\}$ a set of representatives. Let $\tau: H \to \GL{}{V}$ be a linear representation for some finite dimensional vector space $V$ over a field $K$. Define the vector space 
		\begin{equation*}
			W = \bigoplus_{r\in R} rV, \quad rV \coloneqq \{rv \mid v\in V\} \cong V
		\end{equation*}
		Note $rV$ is equipped with the same addition and scalar multiplication as $V$, $rv + rw = r(v+w)$ for all $v,w \in V$ and $trv = r(tv)$ for all $t\in K$. We define the induced representation $\Ind_H^G(\tau): G \to \GL{}{W}$ as the action given by
		\begin{equation*}
			g \cdot \sum_{i=1}^k r_iv_i = \sum_{i=1}^{k} r_i' \tau(h_i)v_i,\quad g r_i = r_i' h_i, \quad g\in G,\, r_i'\in R,\, h_i \in H 
		\end{equation*} 
 	\end{definition}
 	\begin{obs} \label{rmk:induced_dim}
 		It is important to note the dimension of $W$, the vector space where $\Ind_H^G(\tau)$ acts, is $k\dim V$.
 	\end{obs}
 	This notion of induced representations encodes the essence of imprimitivity sets, as such, we have the following reinterpretation of Definition~\ref{def:irred}.
	\begin{definition}\label{def:primitive_alt}
		An irreducible representation $\rho:G \to \GL{n}{\C}$ is \emph{imprimitive} if there exists a proper subgroup $H$ of $G$ and a representation $\tau$ of $H$ such that $\mathrm{Ind}_H^G(\tau)$ is equivalent to $\rho$. If such a representation does not exist, the group is \emph{primitive}.
	\end{definition}
	\begin{obs}
		Following the established notation, when considering a subgroup $H \leq\GL{n}{\C}$ and $\iota: H \hookrightarrow \GL{n}{\C}$ the natural inclusion, we may abuse notation and say $H$ it is reducible, irreducible, imprimitive, etc., instead of $\iota$. In this context, classic literature often refers to $H$ reducible (resp. irreducible) as $H$ intransitive (resp. transitive). Recall that, in general, calling an abstract group $G$ reducible, irreducible, etc. does not make sense for us.
	\end{obs}
	\begin{example}
		We provide examples of imprimitive groups, one monomial and another non-monomial. Let $I_2$ the identity $2\times 2$ matrix and $\omega$ a primitive 5th root of unity. Consider
		\begin{gather*}
			H_1 = \left\langle \begin{pmatrix}
				1 & 0 & 0 \\
				0 & 0 & -1 \\
				0 & 1 & 0
			\end{pmatrix}, 
			\begin{pmatrix}
				0 & 1 & 0 \\
				0 & 0 & 1 \\
				1 & 0 & 0 \\
			\end{pmatrix} \right\rangle, \, h_1 = 
			\begin{pmatrix}
			0 & -1 \\
			1 & 0 \\
			\end{pmatrix}, 
			h_2 = \begin{pmatrix}
			-1 & \omega \\
			\omega^4 &0 \\
			\end{pmatrix}, \\
			H_2 = \left\langle 
			\begin{pmatrix}
			0 & 1 \\
			1 & 0 \\
			\end{pmatrix}\otimes_\Kr I_2, 
			\begin{pmatrix}
				1 & 0 \\
				0 & -1 \\
			\end{pmatrix} \otimes_\Kr I_2, 
			I_2\otimes_\Kr h_1, I_2\otimes_\Kr h_2\right\rangle.
		\end{gather*}
		The group $H_1$ is monomial with imprimitivity set $\{\langle e_1\rangle, \langle e_2\rangle, \langle e_3\rangle\}$, being $e_i$ the canonical basis vectors, and $H_1 \cong A_4\times C_2$. Group $H_2$ is 2-nomial with imprimitivity set $\{\langle e_1, e_2\rangle,\langle e_3,e_4\rangle\}$, $H_2 \cong D_4.A_5$ is \verb|SmallGroup(480,957)| in GAP \cite{GAP4} (in addition, $H_3 = \langle h_1, h_2\rangle$ is a primitive group with $H_3 \cong \SL{2}{5}$).
	\end{example} 
	In many cases, deciding if a linear group is primitive or imprimitive can be difficult. The weaker notion of \emph{quasi-primitivity} is sometimes introduced.
	\begin{definition} \label{def:quasi-primitive}
		An irreducible representation $\rho:G \to \GL{n}{\C}$ is \emph{quasi-primitive} if, for every normal subgroup $N \trianglelefteq G$, the restricted representation $\rho|_N$ is homogeneous, i.e. it is a direct sum of copies of a single irreducible representation (up to equivalence) of $N$.
	\end{definition}
	\begin{obs}
		In our case, we have taken all definitions to require irreducibility on the representation. This is the usual approach to properly separate the study of reducible and irreducible representations, also known as intransitive and transitive in the older literature.
	\end{obs}
	\begin{obs}
		Primitive representations are quasi-primitive due to Clifford's Theorem.
	\end{obs}
	The following result due to \citeauthor{Berger1975} shows that we need only distinguish primitive and quasi-primitive representations when dealing with non-solvable groups.
	\begin{prop}[{\citeauthor{Berger1975} \cite[Main theorem]{Berger1975}}]
		Let $\rho(G) \leq \GL{n}{\C}$ be a finite solvable group. Then $\rho$ is primitive if and only if $\rho$ is quasi-primitive. 
	\end{prop}
	
	So far all definitions have been for linear representations. In order to translate them naturally into the projective setting, we need only observe that they respect projective equivalence.
	\begin{lemma} \label{lemma:proj_equiv}
		Let $G$ be a group and let $\rho_i:G \to \GL{n}{\C}, \, i=1,2$ be linear representations that are projectively equivalent, i.e. if $\pi: \GL{n}{\C} \to \PGL{n}{\C}$ is the projection $\pi \circ \rho_1$ and $\pi\circ\rho_2$ are conjugate by an element of $\PGL{n}{\C}$.
		\begin{enumerate}[(i)]
			\item $\rho_1$ is irreducible if and only if $\rho_2$ is irreducible.
			\item $\rho_1$ is imprimitive if and only if $\rho_2$ is imprimitive.
			\item $\rho_1$ is quasi-primitive if and only if $\rho_2$ is quasi-primitive.
		\end{enumerate}
	\end{lemma}
	\begin{proof}
		Let $f\in\GL{n}{\C}$ and $\mu:G\to \C^\times$ such that $\rho_1(g) = \mu(g)f\rho_2(g)f^{-1}$ for all $g\in G$. Let $V$ be a vector space of $\C^n$, then $\rho(g)(V) = \mu(g)f\rho_2(g)f^{-1}(V) = f\rho_2(g)f^{-1}(V)$. Since irreducibility, primitivity and quasi-primitivity are stable under conjugation (change of basis), so are under projective equivalence.
	\end{proof}
	\begin{definition}
		Let $G$ be a group and $G^*$ a representation group. We say a projective representation $\bar\rho:G\to \PGL{n}{\C}$ satisfies one of the properties in Lemma~\ref{lemma:proj_equiv}(i),(ii),(iii) if a lifting $\rho:G^*\to \GL{n}{\C}$ does. 
	\end{definition}
	\begin{obs}
		Alternatively, we can translate irreducibility, primitivity and quasi-primitivity to the projective setting is to consider the natural action of $\PGL{n}{\C}$ over the vector spaces of $\C^n$, which is the same as the one defined by $\GL{n}{\C}$, except it is not well defined for vectors, only for vector spaces. In the definitions for the linear case, all three properties depend on the action of the group on subspaces of $\C^n$ via their representations. Thus, since lifting from $\PGL{n}{\C}$ to $\GL{n}{\C}$ only involves products by scalars, whose action is trivial on vector spaces, all three definitions \ref{def:irred}, \ref{def:imprimitive} and \ref{def:quasi-primitive} translate naturally to the projective setting and are coherent under liftings. 
	\end{obs}

	In order to classify the finite subgroups of $\GL{n}{\C}$ or $\PGL{n}{\C}$ for a fixed $n$ using \emph{conjugation} as an equivalence relation, the groups are usually separated into the types shown in List \ref{list:types_of_groups}. If item 2.(b) is taken to be ``primitive'', then the types do not overlap. In the case where it is taken to be ``quasi-primitive'', some groups may be quasi-primitive and imprimitive, for instance, $A_6$ has a quasi-primitive projective representation of degree 6 induced by its triple cover $3.A_6$ that is imprimitive (see Table~\ref{tab:pgl_6}).\\
	
	\begin{figure}[h]
			\renewcommand\figurename{List}
		\caption{Types of finite subgroups of $\GL{n}{\C}$ or $\PGL{n}{\C}$}
		\label{list:types_of_groups}
		\begin{enumerate}
			\item Reducible (intransitive).
			\item Irreducible (transitive).
			\begin{enumerate}
				\item Imprimitive.
				\begin{enumerate}
					\item Monomial.
					\item Non-monomial.
				\end{enumerate}
				\item (Quasi-)Primitive.
				\begin{enumerate}
					\item Simple and extensions.
					\item With a reducible normal subgroup and extensions.
					\item With an imprimitive normal subgroup and extensions.
				\end{enumerate}
			\end{enumerate}
		\end{enumerate} 
	\end{figure}
	
	The study of reducible groups can be reduced to studying the irreducible groups of lower degree, as any reducible representation breaks down into irreducible representations of smaller degree. Thus, the classification efforts largely focus on irreducible groups. \\
	
	In the following subsections, we cover some of the main results concerning the structure of imprimitive and primitive groups, including some relatively recent developments. Many of the properties of these groups were first studied for $\GL{n}{\C}$ and then translated into $\PGL{n}{\C}$. As such, most of the stated results will be for $\GL{n}{\C}$.
	
	\subsection{Structure and classification of imprimitive groups} \label{subsec:structure_imprimitive}
	
	The main result characterizing the general structure of imprimitive groups is due to \citeauthor{Blichfeldt1917} \cite[\S 60, Theorem 1]{Blichfeldt1917}.
	\begin{theorem}[{\citeauthor{Blichfeldt1917} \cite[\S 60, Theorem 1]{Blichfeldt1917}}]  \label{thm:imprimitive_structure}
		Let $\rho(G) \leq \GL{n}{\C}$ be an imprimitive group. Let $I(\rho) = \{V_1,\dots,V_k\}$ be a set of imprimitivity of $\rho(G)$. Consider the sets $G_i = \{\rho(g)\in \rho(G)\,|\, \rho(g)V_i = V_i\}, \, i=1,\dots, k$ and $N = \bigcap_{i=1}^kG_i$. Let $\sigma_g : I(\rho) \to I(\rho)$ be the permutation $\sigma_g(V_i) = \rho(g)V_i$. Then $\sigma: G \to S_k$ sending $\sigma(g) = \sigma_g$ is a group morphism, $\sigma(G)$ is a transitive permutation group of degree $k$ and $\ker(\sigma) = N$. In particular, $G_i$ are conjugate in $G$, $[G:G_i]=k$ and $\dim(V_i) = d$ for all $i$ so $n=dk$.
	\end{theorem}
	\begin{obs}
		\citeauthor{Dixon1971} \cite[Theorem 4.2B]{Dixon1971} proves the statement for imprimitive groups of $\GL{}{V}$, $V$ a vector space over any field $F$.
	\end{obs}
	
	The previous result implies that when the dimension $n=p$ is a prime number, all imprimitive groups are monomial. Furthermore, the study of monomial groups reduces to understanding the transitive permutation groups of degree $p$, whose classification is known, and the possible normal subgroups $N$. This structure has been recently exploited to produce exhaustive classifications for the monomial groups of prime degree.
	
	The classification is split between solvable and non-solvable groups. Recall a group is solvable if there exists a subnormal series $\{1\} = G_0 \trianglelefteq G_1 \trianglelefteq \cdots \trianglelefteq G_k = G$ such that $G_j/G_{j-1}$ is an abelian group for $j=1,\dots, k$.
	
	\begin{classification}[Monomial groups of prime degree] \label{class:monomial}
		Let $G$ be a group, $p$ be a prime number and $K$ an algebraically closed field.
		\begin{enumerate}[(i)]
			\item Let $G$ be non-solvable and $\rho:G\to \SL{p}{K}$ a monomial representation, then $\rho(G)$ has been classified (not listed) by \citeauthor{DixonZalesskii2004Imprimitive} in \cite{DixonZalesskii2004Imprimitive} up to conjugation. (See Remark \ref{rmk:non_solvable})
			\item Let $G$ be solvable and $\rho:G\to \GL{p}{\C}$ a monomial representation, then $\rho(G)$ has been listed by \citeauthor{Flannery2021} in \cite{Flannery2021} up to conjugation. (See Remark \ref{rmk:solvable})
		\end{enumerate}
	\end{classification}
	\begin{obs} \label{rmk:non_solvable}
		In the non-solvable case, \citeauthor{DixonZalesskii2004Imprimitive} in \cite{DixonZalesskii2004Imprimitive} reduce the problem to studying $\SL{p}{\C}$ rather than $\GL{p}{\C}$. The classification is complete for $\GL{p}{\C}$ except for the explicit extensions of the groups to $\GL{p}{\C}$. This restriction does not affect the induced classification of $\PGL{p}{\C}$, for $\PSL{p}{\C} = \PGL{p}{\C}$, see Remark~\ref{rmk:unimodular_lift}.
	\end{obs}
	
	\begin{obs} \label{rmk:solvable}
		In the solvable case, \citeauthor{Flannery2021} in \cite{Flannery2021} provide a fully explicit description of the monomial groups of prime degree, in the sense that representations of the group can be obtained. This description has been incorporated in the computer algebra system Magma \cite{Magma} and is publicly available to any user, see \cite{Monomial2021}. We stress the relevancy of these results, for they are the ``ideal'' form a classification of the finite irreducible subgroups of $\GL{n}{\C}$, in the sense that the explicit representations of the groups can be obtained for any of the groups using \cite{Monomial2021}. This is crucial for potential applications to other fields. Unfortunately, in the non-solvable case such a detailed description is not available, as far as we know.
	\end{obs}
	
	\begin{obs}
		 \citeauthor{Flannery2021} in \cite{Flannery2021} also tackle part of the non-solvable case and describe it except for a particular type of primes, see \cite[\S 10]{Flannery2021} for more details. Their description of the non-solvable case concludes the explicit classification of the finite monomial subgroups of $\GL{p}{\C}$ for $p\leq 11$ and $p=23$. Furthermore, they also provide all irreducible representations of $\GL{2}{\C}$ and $\GL{3}{\C}$, not only the imprimitive ones, as part of their computational implementation in \cite{Monomial2021}.
	\end{obs}
	
	Classification \ref{class:monomial} concludes most of the work on the monomial case for prime dimension. The theory of monomial groups of composite degree and non-monomial groups is far less developed. Some relevant results include the monomial groups of $\GL{4}{\C}$ by \citeauthor{Flannery1999} in \cite{Flannery1999} and the description of the imprimitive non-monomial groups of degree 4 by \citeauthor{Hofling2001} in \cite{Hofling2001}. Unfortunately, as far as the authors know, these classifications are not currently implemented on any computer system.
	
	\begin{obs} \label{rmk:problem_translation}
		Fully translating Classification \ref{class:monomial} into a classification for $\PGL{n}{\C}$ is no easy task \emph{a priori}. If $G$ is a group and $\rho(G) \leq \GL{p}{\C}$ is monomial, it is necessary to identify the subgroup of $Z(G)$ corresponding to $\rho^{-1}(\rho(G)\cap Z(\GL{p}{\C}))$. By the decomposition of Theorem \ref{thm:imprimitive_structure}, this group is in $N = \ker(\sigma)$, the normal subgroup of diagonal matrices of $G$. One would need to asses for each presented family of monomial groups in \cite{Flannery2021} and \cite{DixonZalesskii2004Imprimitive}, what part of the associated group of diagonal transformations is scalar and needs to be quotiented out. Moreover, even if this process is carried out, the serious problem comes with redundancy removal, as projective equivalence is fundamentally different from linear equivalence. As far as the authors know, tasks to asses redundancies in this setting have not been undertaken.
	\end{obs}
	\begin{example}
		Let us exemplify the practical consequences of the problem in Remark~\ref{rmk:problem_translation}. Take the computational implementation of the classification of the solvable monomial groups in \cite{Monomial2021}. Using the program, one can obtain all imprimitive linear groups of degree $p$ and order $N$ both fixed. However, if we wanted all imprimitive projective groups of dimension $p$ and order $N$ both fixed, it is not enough to only search for the linear groups satisfying these conditions. There are mainly two reasons for this. First, the linear groups obtained may contain scalar multiples in their center, so when projecting the image is no longer of size $N$. Second, in the same way, there may be linear imprimitive groups of order greater than $N$ which project to a group of order $N$, and these need to be taken into account.
	\end{example}
		
	As a final remark on the structure of these groups, a theory of invariant polygons for monomial groups of $\GL{n}{\C}$ has been developed by \citeallauthors{KangZhang2009} in \cite{KangZhang2009}. The concept of invariant polygon, first introduced by \citeauthor{Blichfeldt1917}, essentially refers to a collection of $n$ vectors which are permuted by the action of a monomial group (modulo scalars). It can be seen that, for a monomial group $\rho(G)$, there are finitely many $\rho(G)$-invariant polygons. This theory is relevant to Section \ref{sec:primitive_prime}, as it is useful to describe primitive groups with monomial normal subgroups.
	
	\subsection{Structure of (quasi-)primitive groups} \label{subsec:structure_primitive}
	We start the section presenting the first fundamental structure theorem involving quasi-primitive linear groups, originally due to \citeauthor{Blichfeldt1917}.
	
	
	\begin{theorem}[{\citeauthor{Blichfeldt1917} \cite[\S 61, Lemma]{Blichfeldt1917}}] \label{thm:cyclic_center} \label{thm:primitive_center}
		Let $\rho(G) \leq \GL{n}{\C}$ be a (quasi-)primitive linear group. Let $A \leq G$ be a normal abelian subgroup, then $A\leq Z(G)$. Furthermore, $Z(G)$ is cyclic and $\rho(Z(G))$ are scalar matrices.
	\end{theorem}
	\begin{obs}
		\citeauthor{Blichfeldt1917}'s proof is from before quasi-primitive groups were introduced. Nonetheless, his proof also holds for quasi-primitive groups.
	\end{obs}
	\begin{obs} \label{rmk:aschbacher_general}
		Theorem~\ref{thm:cyclic_center} was greatly generalized by \citeauthor{Aschbacher2000} in \cite[\S 1.6]{Aschbacher2000}. In particular, the statement holds for quasi-primitive groups over any field $K$ so long as the representation $\rho:G \to \GL{n}{K}$ is completely reducible (also known as semisimple) and for each $H\trianglelefteq G$ normal subgroup, each irreducible $K[H]$-submodule of $K^n$ is absolutely irreducible. \citeauthor{Aschbacher2000} calls these conditions \emph{AI}. Over an algebraically closed field, this reduces to simply requiring $\rho$ to be completely reducible.
	\end{obs}
	\begin{obs}
		Note that primitive projective groups are not necessarily without normal abelian subgroups, for taking the quotient by the center does not preserve this property. However, it is true that they are centerless.
	\end{obs}

 	\begin{obs} \label{rmk:CR-radical}
 		Groups $G$ without normal abelian subgroups are sometimes called \emph{semisimple} \cite[\S 3.3]{Robinson1995BOOK} or Fitting-free groups. These groups contain a normal subgroup $R\trianglelefteq G$ called the \emph{centerless completely irreducible radical}, which is a direct product of non-abelian simple groups. Furthermore, $\mathrm{Inn}(R)\leq G \leq\mathrm{Aut}(R)$ \cite[3.3.18]{Robinson1995BOOK} and the possible groups fitting between $\mathrm{Inn}(R)$ and $\mathrm{Aut}(R)$ are in bijective correspondence to the conjugacy classes of $\mathrm{Out}(R)$ \cite[3.3.20]{Robinson1995BOOK}. These groups are tightly related to quasi-primitive groups as will be seen in Lemma~\ref{lemma:DZ_cases}(ii) or Theorem~\ref{thm:composite_cases}(2),(4).
 	\end{obs}
 	
 	By Propositions \ref{prop:lift_faithful} and \ref{prop:proj_faithful}, studying the (quasi-)primitive projective representations of the simple groups is the same as studying the (quasi-)primitive linear representations of quasisimple groups with cyclic center. Recall a group $G$ is quasisimple if it is perfect (i.e. $G'=G$) and it is a central extension of a simple group. 
 	
 	The \emph{Classification of Finite Simple Groups} (CFSG) provides valuable information for our problem. The Schur coverings of simple groups are unique and have been extensively studied as part of the CFSG. The ATLAS \cite{ATLAS} contains most of the information about the coverings of the sporadic groups and some low order Lie type groups. All Schur multipliers of the simple groups are known. Furthermore, descriptions and lower bounds for the degrees of the projective representations of groups of Lie type were given by Landazuri, \citeauthor{Seitz1993} in \cite{Seitz1974} \cite{Seitz1993}.
 	
 	\begin{example}
 		We provide an example of a lower bound for the degree of a projective re\-pre\-sentation given in \cite{Seitz1993}. Let $G(q)$ be a finite simple group of Lie type over a field of order $q$. We denote by $l(G(q))$ the smallest integer $d>1$ such that $G(q)$ has a projective representation over a field of characteristic coprime to $q$. Let $n\geq 3$, then $l(\PSL{}{n,q}) = \frac{q^n-1}{q-1}-n$ with the following exceptions $l(\PSL{}{3,2}) = 2$, $l(\PSL{}{3,4}) = 4$, $l(\PSL{}{4,2}) = 7$ and $l(\PSL{}{4,3}) = 26$.
 	\end{example}
 	
 	A classification of the complex irreducible representations of the quasisimple groups for re\-la\-ti\-ve\-ly small degree due to \citeauthor{Tiep1996} \cite{Tiep1996} was surveyed, alongside many other results involving the degrees of representations of these groups, by the same authors in \cite{Tiep2000}. In particular, see \cite[Theorem 6.1]{Tiep2000} or Theorem~\ref{teoTZ96} for the classification.
 	
 	Nevertheless, great effort has been made to prove structural results involving $p$-Sylow subgroups for (quasi)-primitive groups. We present some of them alongside relevant properties of (quasi-)primitive groups derived from the former.
 	
 	\begin{theorem}[\citeauthor{Blichfeldt1917}, \citeauthor{Brauer1967}] \label{thm:bound_primes}
 		Let $\rho(G^*) \leq \GL{n}{\C}$ be a finite linear group with order $g=|G^*|$. Then for each prime $p\mid g$ with $p> 2n+1$, $G^*$ has a unique $p$-Sylow subgroup. In particular, if $G^*$ is a covering of $G$ and $\rho$ is (quasi-)primitive, then there does not exist $p>2n+1$ such that $p\mid [G^*:Z(G^*)]$.
 	\end{theorem}
 	\begin{obs}
 		\citeauthor{Blichfeldt1917}'s original result \cite[\S 64, Theorem 5]{Blichfeldt1917} gave a bound $p>(2n+1)(n-1)$. \citeauthor{Brauer1967} \cite[2C]{Brauer1967} deduced the uniqueness of the relevant $p$-Sylow subgroup for $p>2n+1$ reducing the bound. A modernized proof of \citeauthor{Blichfeldt1917}'s result was given by \citeauthor{Dixon1971} in \cite[Theorem 5.5]{Dixon1971}.
 	\end{obs}
 	\begin{obs} \label{rmk:proj_to_lin_by_quot_center}
 		In the sequel, many results concerning (quasi)-primitive groups $\rho(G)$ will be stated for $G/Z(G)$ instead of $G$. To understand the implications of this, recall $\rho(Z(G))$ are scalar matrices by Theorem~\ref{thm:primitive_center}. Thus, if $\pi:\GL{n}{\C} \to \PGL{n}{\C}$ is the canonical projection, $\pi \circ \rho: G \to \PGL{n}{\C}$ has as its image $G/Z(G) \cong (\pi \circ \rho)(G)$. Thus, it defines a faithful projective representation $\bar{\rho}: G/Z(G) \to \PGL{n}{\C}$ that is lifted by $\rho: G\to \GL{n}{\C}$, meaning $G$ is in fact a covering of $G/Z(G)$ by $Z(G)$. If we call $\bar{G} = G/Z(G)$, then results given for $\bar{G}$ concern the projective primitive group while $G$ refers to coverings of $\bar{G}$ by $Z(G)$. In our notation, we could write $G = Z(G).\bar{G}$.
 	\end{obs}
 	
 	Further results by \citeauthor{Brauer1967} \cite[2D]{Brauer1967} limit the presence of primes larger than the degree.
 	\begin{theorem}[\citeauthor{Brauer1967} \cite{Brauer1967}] \label{thm:bound_primes_ged_dim}
 		Let $\rho(G) \leq \GL{n}{\C}$ be a primitive group and $p>n+1$ a prime. Then $p$ is the only prime larger than $n+1$ dividing $|G|$, $p^2\nmid |G|$ and if $p=2n+1$ then $G/Z(G) \cong \PSL{}{2,p}$.
 	\end{theorem}
 	
 	There are many results concerning the structure of the $p$-Sylow groups for specific degrees, many obtained after the original work of \citeauthor{Blichfeldt1917}. We list some of them and apologize for the many omissions.
 	
 	\begin{prop}[\citeauthor{Bra42}, \cite{Bra42}] Let $\rho(G) \leq \mathrm{GL}_{(p-1)/2}(\C)$ be a primitive group for some prime $p$. Then, either $G/Z(G)\cong \mathrm{PSL}(2,{p})$ or $G$ has an abelian normal $p$-Sylow subgroup.
 	\end{prop}
 	 
 	\begin{prop}[\citeauthor{Hay63}, \cite{Hay63}] Let $\rho(G) \leq \mathrm{GL}_{(p+1)/2}(\C)$ be a primitive group for some prime $p\geq 7$. Then either $G$ has an abelian normal $p$-Sylow subgroup or else $G/Z(G)\cong \mathrm{PSL}(2,{p})$.
 	\end{prop}

 	\begin{prop}[Feit, \cite{Feit64, Feit67}] Let $\rho(G) \leq \mathrm{GL}_{p-2}(\C)$ be a finite irreducible group with $p$ prime, then, either $G$ has an abelian normal $p$-Sylow subgroup or $p$ is a Fermat prime and $G/Z(G)\cong \SL{}{{2},p-1}$.
 	\end{prop}
 	 
 	Exploiting the structure of the $p$-Sylow groups and their intersections, \citeauthor{Blichfeldt1917} was able to prove a bound for the maximum power of certain primes in a primitive group.
 	 
 	\begin{theorem}[\citeauthor{Blichfeldt1917}, \cite{Blichfeldt11}] \label{thm:bound_primes_order}
 		Let $\rho(G) \leq \SL{n}{\C}$ be a primitive group, let $p$ be a prime number such that $p\nmid n$. Consider $k\geq 1$ such that $p^k \mid |G|$, then $p^{k} \mid n!p^{n-1}$.
 	\end{theorem}
 	
 	Following a different approach, \citeauthor{Blichfeldt1917} was able to provide a generic bound for any prime. This alternative approach shifts the focus towards the eigenvalues of a representation, it is developed in \cite[\S\S 69-74]{Blichfeldt1917}. As an example, take the following result.
 	
 	\begin{lemma}[\citeauthor{Blichfeldt11}, \cite{Blichfeldt1917}]
 		Let $\rho(G) \leq GL_n(\mathbb{C})$ be a finite quasi-primitive group and let $g\in G$ be such that all the eigenvalues of $\rho(g)$ lie within $\pi/3$ of one of them on the unit circle, then $g\in Z(G)$.
 	\end{lemma}
 	
 	These ideas culminate in the following theorem.
 	\begin{theorem}[{\citeauthor{Blichfeldt1917}, \cite[\S 74]{Blichfeldt1917}}] \label{thm:prime_bounds_general}
 		Let $\rho(G) \leq GL_n(\mathbb{C})$ be a finite quasi-primitive group, let $p$ be a prime number and $k\geq 1$ such that $p^k\mid |G|$. Then $p^k \leq (n!)_p6^{n-1}$ where $(n!)_p$ denotes the maximum power of $p$ dividing $n!$.
 	\end{theorem}
 	\begin{obs}
 		In \cite{Blichfeldt1917}, Theorem~\ref{thm:prime_bounds_general} is stated with a $5$ instead of a $6$. As \citeauthor{Bra42} notes in \cite[p.74]{Brauer1967}, this reduction from $5$ to $6$ is presented as pendent of publication. It seems, however, that it never ended up being published. 
 	\end{obs}
 	
	The previous theorems involve bounds on the prime numbers dividing the order of a primitive group. Recent results by \citeauthor{Collins2008} provide sharp bounds for the order of the primitive groups themselves.
	
	\begin{theorem}[{\citeauthor{Collins2008} \cite[Theorem A]{Collins2008}}] \label{thm:collins_bounds}
		Let $\rho(G) \leq \GL{n}{\C}$ be a primitive group, then $[G:Z(G)]\leq (n+1)!$ if $n>12$ or $n=10,11$. In particular, the bound is achieved when $G' \cong A_{n+1}$ and $G/Z(G)\cong S_{n+1}$.
	\end{theorem}
	\begin{obs}
		The exceptions to the bound for $n\leq 12$ in Theorem~\ref{thm:collins_bounds} are also provided in a table in \cite[Theorem A]{Collins2008}.
	\end{obs}
	
	\begin{obs}\label{rmk:translation}
		All results stated so far in this subsection trivially generalize to statements over $\PGL{n}{\C}$ through Theorem~\ref{thm:primitive_center} and Remark~\ref{rmk:proj_to_lin_by_quot_center}, as they are already stated quotient the center of the primitive group.
	\end{obs}
	
	\begin{obs}
		Theorem~\ref{thm:collins_bounds} alongside Theorem~\ref{thm:cyclic_center} imply that there are finitely many finite primitive linear groups on $\PGL{n}{\C}$ up to conjugation, for the order of the groups is upper bounded, and there are finitely many groups of finite fixed order. This also extends to primitive groups over $\SL{n}{\C}$, as $Z(\SL{n}{\C}) \cong C_n$. In contrast, for imprimitive groups, the normal abelian diagonal subgroups (given by the kernel of $\sigma$ described in Theorem~\ref{thm:imprimitive_structure}) may grow arbitrarily large. In this case we talk about families of imprimitive groups, see \cite{Flannery2021} and Classification~\ref{class:monomial} for examples.
	\end{obs} 
	
	\citeauthor{Collins2008} uses Theorem~\ref{thm:collins_bounds} to provide sharp bounds for any finite linear complex group in \cite[Theorem B]{Collins2007}.
	
	\begin{theorem}[\citeauthor{Collins2007} \cite{Collins2007}]
		Let $\rho(G) \leq \GL{n}{\C}$ be a linear group and $H\trianglelefteq G$ be an abelian normal subgroup with maximal order. Then $[G:H] \leq f(n)$ where
		\begin{enumerate}[(i)]
			\item $f(n) = (n+1)!$ if $n\geq 71$ or $n= 63,65,67,69$. The bound is realized when $G' \cong A_{n+1}$ and $G/Z(G)\cong S_{n+1}$.
			\item $f(n) = 60^rr!$ for $n=2r$ or $n=2r+1$ if $20\leq n\leq 70$ and $n\neq 63,65,67,69$. The bound is realized when $G$ has a normal subgroup $H$ such that $H = Z(H).E(H)$ is a central extension and $E(H)$ is a direct product of $r$ copies of $\SL{2}{5}$ and $G/H\cong S_r$.
		\end{enumerate}
	\end{theorem}
	For $n\leq 20$, \cite[Theorem D]{Collins2007} specifies groups whose indices provide the values of $f(n)$. These result concludes a program started by Jordan on 1878 to identify the function $f(n)$ bounding the indices $[G:H]$ for finite subgroups of $\GL{n}{\C}$.

	\newpage
	\section{Primitive subgroups of $\SL{n}{\C}$ and $\PGL{n}{\C}$} \label{sec:primitive_groups}
	
	In this section we focus on primitive and quasi-primitive groups and what is known of their classifications, both for $\PGL{n}{\C}$ and $\SL{n}{\C}$ for some fixed $n \geq 2$. The first subsection focuses on the prime dimensional case, where most effort has been done. Throughout the second subsection we discuss some of the results known when $n$ is composite.
	
	\subsection{The prime dimensional case} \label{sec:primitive_prime}
	
	Let $p$ be a prime, in this subsection we focus on the finite primitive subgroups of $\PGL{p}{\C}$ and $\SL{p}{\C}$. We present two main approaches to the classification. The classic approach, which is pre-CFSG, and the post-CFSG approach, where we recall CFSG means classification of finite simple groups. 
	
	The classic approach is centered on a program, started by \citeauthor{Blichfeldt1917}, to describe the finite primitive groups of $\PGL{n}{\C}$ for small $n$. This program started around 1917 with \citeauthor{Blichfeldt1917}'s classification of $n=3,4$ and ended with \citeauthor{Rob81}'s PhD thesis description of $n=11$ in 1981. In particular, this program was developed before the CFSG and classified the primitive groups for $2\leq n\leq 11$. A description of this classification can be found in Tables~\ref{tab:pgl_2}-\ref{tab:pgl_7} in Section~\ref{sec:tables} and in Classification~\ref{class:primitive_small_p}.
	
	The post-CFSG approach was led by \citeauthor{DixonZalesskii1998Primitive} among many other mathematicians. The idea is to exploit all the information about the projective representations of finite simple groups and quasi-simple groups to describe the primitive groups with non-abelian socle. Recall that the socle of a group is the subgroup generated by all the minimal normal subgroups. The abelian socle case has been explicitly developed by \citeallauthors{KangZhang2009} in \cite{KangZhang2009} to provide a construction for explicit representations of such groups via the development of a theory of invariant polygons, first used by \cite{Blichfeldt1917}.
	
	\subsubsection{The classic approach} \label{subsubsec:classic_approach}
	When considering finite subgroups $G$ of $\GL{p}{\C}$, we always have $p$ divides the order of $G$. This is a standard result in representation theory. This always enables us to consider a non-trivial $p$-Sylow subgroup of $G$. We start studying these special subgroups.
	
	\begin{lemma}[Sibley \cite{Sib74}]\label{BraLinSib}
		Let $\rho(G)\leq \SL{p}{\C}$ be a finite quasi-primitive group and $p\geq 5$ a prime. Take $P$ a Sylow $p$-group of $G$. Then,
		\begin{enumerate}[(i)]
			\item if $|P|=p^2$, then $G\cong G_1\times C_p$ where $G_1$ is a group whose $p$-Sylow has order $p$,
			\item if $|P|=p^3$, then $P$ is normal in $G$ and $G/P$ is isomorphic to a subgroup of $\SL{}{2,p}$,
			\item if $|P|=p^4$, then $P$ has a subgroup $Q$ of index $p$ which is normal in $G$ and $G/Q$ is isomorphic to a subgroup of $\SL{}{2,p}$,
			\item if $|P|\geq p^5$, no such $G$ does not exists.
		\end{enumerate}
	\end{lemma}
	\begin{obs}
		Lemma~\ref{BraLinSib} appears as well-known in the introduction of \cite{Sib74} (except $|P|=p^3$), and refers to \citeauthor{Bra42} in \cite{Brauer1967} for $|P|=p^2$ and the work of \citeauthor{Lin70} \cite{Lin71}\cite{Lin70} for $|P|\geq p^4$. It also claims that Feit also obtained independently such results (but they are unpublished). The case $|P|=p^3$ is the contents of the paper \cite{Sib74}.
	\end{obs}
	\begin{obs}
		Note that cases (ii) and (iii) imply that $P$ must be non-abelian. Otherwise, $G$ could not have a quasi-primitive representation by Theorem~\ref{thm:primitive_center}. Note that over $\SL{p}{\C}$ we have $|Z(\rho(G))|$ is either $1$ or $p$, so $P$ does not fit in $Z(G)$ for any of these cases.
	\end{obs}
	We are missing a description of $G$ when $|P|\leq p$. The following lemma due to Brauer covers this case.
	\begin{lemma}[Brauer]\label{lemma:p_sylow_p} Let $G$ be a finite primitive of prime degree $p$ with $p\geq 3$, with an abelian $p$-Sylow subgroup $P$, then $G'$ is a non-abelian simple group. Furthermore, if $|P| = p$ then $Z(G) = 1$.
	\end{lemma}
	\begin{obs}
		We lack the original reference of Brauer, but Robinson gave an alternative proof in \cite[Lemma 2.1.2]{Rob81}. The implication: if $|P| = p$ then $Z(G) = 1$; can be found in \cite[4a]{Brauer1967}.
	\end{obs}
	
	From Lemma~\ref{BraLinSib} and Lemma~\ref{lemma:p_sylow_p}, when $p\geq 5$, there are two possible structures for a primitive group $\rho(G)$ with $\rho:G\hookrightarrow \SL{p}{\C}$:
	\begin{enumerate}[(A)]
		\item \label{case:A} $G$ has a normal subgroup $N$ of order $p^3$ and $G/N$ is isomorphic to a subgroup of $\SL{}{2,p}$.
		\item \label{case:B} $G'$ is a non-abelian simple group.
	\end{enumerate}
	\begin{obs} \label{rmk:case_i}
		Case \ref{case:A} can be further specified to show that the normal subgroup $N$ of $G$ with order $p^3$ is in fact an extra special group such that $|Z(N)|=p$ and $N/Z(N)$ is an elementary abelian group \cite[\S 3, Proof Lemma 1.1]{DixonZalesskii1998Primitive}. 
	\end{obs}
	
	\let\originalthetheorem\thetheorem
	\renewcommand\thetheorem{\thesection.\number\numexpr\value{theorem}bis}
	\addtocounter{theorem}{-1}
	\begin{obs} \label{rmk:bis}
		Case \ref{case:B} corresponds to the groups with abelian $p$-Sylow subgroups by Lemma~\ref{lemma:p_sylow_p}. Furthermore, we can reduce to studying $G$ with $|P|=p$ by Lemma~\ref{BraLinSib} and consider $Z(G) = 1$. This case is closely related to the discussion in Remark~\ref{rmk:CR-radical}, where $R = G'$ and so the relevant groups in this case are the groups $G$ such that $G'\leq G\leq\mathrm{Aut}(G')$. One can also view this from the perspective of group extensions of $G'$ by $\mathrm{Out}(G')$. Furthermore, this also shows that projective representations of groups in case \ref{case:B} are all given by trivial coverings, i.e. obtaining the list of groups in case \ref{case:B} reduces to finding all non-abelian simple groups $H$ and all groups $G$ such that $H\leq G \leq \mathrm{Aut}(H)$ and determining which of these $G$ are primitive.
	\end{obs}
	\let\thetheorem\originalthetheorem
	
	By Remark\ref{rmk:bis}, the efforts in classifying primitive groups of degree $p$ when the $p$-Sylow subgroup is abelian, i.e. case \ref{case:B}, is centered in obtaining the classification for simple groups $G$ with an abelian $p$-Sylow subgroup. For small primes $2\leq p\leq 11$, the description of case \ref{case:B} is as follows.
	\begin{theorem} \label{thm:groups_small_p}
		Let $\rho(G)$ be a finite primitive group of $SL_n(\C)$ with $G'$ non-abelian simple, then
		\begin{enumerate}
			\item for $p=2$, $|Z(G)|=2$ and $G/Z(G)$ is isomorphic to $A_5$. Otherwise,  $|Z(G)|=2$ and $G/Z(G)$ is isomorphic to $A_4$ or $S_4$.
			\item for $p=3$, $G/Z(G)$ is isomorphic to $A_5$, $A_6$ or $\mathrm{PSL}(2,7)$.
			\item (Brauer \cite{Brauer1967}) for $p=5$, $G/Z(G)$ is isomorphic to $S_5$, $A_6$, $S_6$, $\mathrm{PSL}(2,11)$, $\mathrm{PSU}(4,2)$.
			\item (Wales \cite{Wales1969},\cite{Wales1970}) for $p=7$, $G/Z(G)$ is isomorphic to $A_8$,$S_8$, $\mathrm{PSL}(2,13)$, $\mathrm{PSp}(6,2)$, $\mathrm{PGL}(2,7)$, $\mathrm{PSL}(2,8)$, $\mathrm{R}(3)$, $\mathrm{PSU}(3,3)$ or $\mathrm{G}_2(2)$.
			\item (Robinson \cite{Rob81}) for $p=11$, $G'\cap Z(G)=1$, $[G:G'\times Z(G)]\leq 2$ and $G'$ is isomorphic to one of the following groups: $A_{12}$, $M_{12}$, $\mathrm{PSL}(2,11)$, $\mathrm{PSL}(2,23)$, or $\mathrm{PSU}(5,2)$. 
		\end{enumerate}
	\end{theorem}
	\begin{corrections} \label{correction:small_p_1}
		For $p=5$, Brauer also considered $A_5$ as a primitive group. However, \citeauthor{DixonZalesskii1998Primitive}'s classification (see Theorem~\ref{thm:socle_arb_p}) states that $A_5$ appears as an imprimitive group. A similar problem occurs for $p=7$, where \citeauthor{Wales1969} considered $\PSL{}{2,7}$ a primitive group, but \citeauthor{DixonZalesskii1998Primitive}'s classification states it appears as an imprimitive group. To clarify these contradictions, we provide a short proof that $A_5$ is not primitive.
	\end{corrections}
	\begin{lemma} \label{lemma:A_5_imprimitive}
		The unique degree $5$ complex irreducible representation of $A_5$ is imprimitive.
	\end{lemma}
	\begin{proof}
		We consider the characterization of imprimitivity through induced representations, see Definition~\ref{def:primitive_alt}. 
		
		We want to induce the unique degree $5$ representation of $A_5$ using a proper subgroup. From Remark~\ref{rmk:induced_dim}, we see that it must be induced by a subgroup $H \leq A_5$ with $[G:H] = 5$ and through a linear character, meaning a 1-dimensional representation of $H$. 
		
		With this information, we easily see that $H \cong A_4$ is the only option. Furthermore, note that $A_4$ has 2 non-equivalent linear characters. Let $\tau: A_4 \to \C$ be one of the two non-trivial linear
		characters and consider $\chi = \Ind_H^G(\tau)$. 
		
		We now proceed to analyze how $\chi$ decomposes into a direct sum of irreducible representations, we call these the \emph{constituents} of $\chi$. Note that $A_5$ is simple non-abelian, so it has no non-trivial linear characters, meaning the trivial representation is the only possible degree 1 constituent of $\chi$.
		
		Recall that Frobenius reciprocity states, for some class function $\varphi$ of $G$:
		\begin{equation*}
			\langle \Ind_H^G(\tau), \varphi\rangle_G = \langle\tau, \Res_H^G(\varphi)\rangle_H
		\end{equation*}
		Choosing $\varphi = 1_G$, the trivial character of $G$, we see
		\begin{equation*}
			\langle \chi, 1_G\rangle_G = \langle\tau, 1_H\rangle_H = 0
		\end{equation*}
		so the trivial character is not a constituent of $\Ind_H^G(\tau)$. Since the degrees of the irreducible representations of $A_5$ are $1,3,3,4,5$. In particular, there is no degree 2 irreducible representation, so the only option is for $\chi = \Ind_H^G(\tau)$ to be irreducible, and so $A_5$'s degree 5 representation is imprimitive. 
	\end{proof}
	\begin{obs}
		A similar argument can be used to show $\PSL{2}{7}$ is imprimitive, using the subgroup $[\PSL{2}{7}: S_4] = 7$. In fact, the simple group $M_{11}$ also famously has an imprimitive irreducible representation of degree $11$, whose monomial form is known (see the GAP package \cite{AtlasRep2.1.6} or \cite{ATLASReps}). The same procedure as in Lemma~\ref{lemma:A_5_imprimitive} shows $M_{11}$ is imprimitive using its index $11$ subgroup $M_{10}$.
	\end{obs}
	\begin{obs}
		We can also take this opportunity to make use of the recent result of \citeauthor{Flannery2021} \cite{Flannery2021} (see also Classification~\ref{class:monomial}) to show that $A_5$ and $\PSL{2}{7}$ are imprimitive.
		Before proceeding to do this, we believe it is relevant to provide more context on Correction~\ref{correction:small_p_1} via the statement of Theorem~\ref{thm:socle_arb_p}. Thus, we continue this discussion in Correction~\ref{correction:small_p_2}.
	\end{obs}
	
	An explicit list of the groups conforming case \ref{case:A} was not given in the cited classifications \cite{Brauer1967}, \cite{Wales1969}, \cite{Wales1970}, \cite{Rob81}, only the description provided in Remark~\ref{rmk:case_i} is given. However, \citeallauthors{KangZhang2009} in \cite{KangZhang2009} developed a theory to obtain the representations of the groups of case \ref{case:A}, we discuss it at the end of this section (after Correction~\ref{correction:small_p_2}). For now, we can sum up the presented results in the following statement.
	
	\begin{classification}[Primitive groups of small prime degree]\label{class:primitive_small_p}
		Let $G$ be a finite group. Let $p$ be a prime number such that $2\leq p\leq 11$. Let $\rho: G \to \SL{p}{\C}$ be a primitive linear representation of $G$.
		\begin{enumerate}[(i)]
			\item If $G$ has an abelian $p$-Sylow subgroup, then $G'$ is simple and $G$ is one of the groups in Theorem~\ref{thm:groups_small_p}. For $p=11$, the description of the groups is not explicit and we may use Remark~\ref{rmk:bis} to obtain all groups.
			\item If $G$ has a non-abelian $p$-Sylow subgroup $P$, then $|P|\geq p^3$. For $p=2$ there are none, for $p=3$ they were described by \citeauthor{Blichfeldt11} in \cite[\S 79]{Blichfeldt1917}. For $p=5,7$ these were described in \cite[Theorems A.3, A.6]{KangZhang2009}.
		\end{enumerate}
	\end{classification}
	\begin{obs}
		The given classification over $\SL{p}{\C}$ projects into a classification of $\PGL{p}{\C}$ for the groups in Classification~\ref{class:primitive_small_p}(i), as Theorem~\ref{thm:groups_small_p} already gives the groups quotient the center, so these are all the groups by Theorem~\ref{thm:cyclic_center} and Remark~\ref{rmk:proj_to_lin_by_quot_center}. For Classification~\ref{class:primitive_small_p}(ii), one needs to quotient out the center to obtain the groups in $\PGL{p}{\C}$. An explicit list of these groups for $2\leq p \leq 7$, alongside additional information, can be found in Section~\ref{subsec:tables} tables \ref{tab:pgl_2},\ref{tab:pgl_3},\ref{tab:pgl_5} and \ref{tab:pgl_7}.  
	\end{obs}
	
	As a historical note, the first tables collecting all such results covered in Classification~\ref{class:primitive_small_p} for $2\leq p\leq 7$ were given by Feit \cite{Feit1971} in 1971. The tables we provide in Section~\ref{sec:tables} are a refinement of his tables, with modern notation and the additional results of \cite{KangZhang2009} covering the ``$I_p$ groups'', as Feit called them in \cite[p.76]{Feit1971}, which refer to those in Classification~\ref{class:primitive_small_p}(ii).
	
	\subsubsection{The post-CFSG approach} \label{subsubsec:post_approach}
	In the sequel, we describe the approach of \citeauthor{DixonZalesskii1998Primitive} to the classification of the primitive groups of $\SL{p}{\C}$ with non-abelian socle. Note that all presented results immediately translate into a classification of $\PGL{p}{\C}$ by Remarks~\ref{rmk:proj_to_lin_by_quot_center} and \ref{rmk:translation} as we mostly consider groups quotient their center. The abelian socle case is also treated later in the section. 
	
	This approach aims to explicitly describe the socles of these groups for each $p$. We denote the socle of a group $G$ as $\mathrm{soc}(G)$. Recall it is the subgroup generated by all the minimal normal subgroups of $G$. In the context of the primitive groups, $\mathrm{soc}(G/Z(G))$ coincides with the centerless completely irreducible (CR) radical of $G$ described in Remark~\ref{rmk:CR-radical}. 
	
	We start reformulating the two cases for $G$ stated after Lemma~\ref{lemma:p_sylow_p} in a more suitable way for this context.
	\begin{lemma}[{\citeauthor{DixonZalesskii1998Primitive} \cite[Lemma 1.1]{DixonZalesskii1998Primitive}}] \label{lemma:DZ_cases}
		Let $\rho(G)\leq \SL{p}{\C}$ be a finite primitive group and write $S \coloneqq \mathrm{soc}(G/Z(G))$. Then one of the following situations holds:
		\begin{enumerate}[(i)]
			\item $S$ is an elementary abelian $p$-group of order $p^2$ and $G/Z(G)$ is isomorphic to a subgroup $H$ of an extension of $S$ by $\SL{}{2,p}$ (split if $p>2$),
			\item $S$ is a non-abelian simple group and $G/Z(G)$ is isomorphic to a subgroup of $\mathrm{Aut}(S)$.
		\end{enumerate}
	\end{lemma}
	
	The possible non-abelian simple groups $S$ appearing in case (ii) were described in terms of $p$ by \citeauthor{DixonZalesskii1998Primitive} in \cite{DixonZalesskii1998Primitive} and \cite{DixonZalesskii2008Corrigendum}. The idea of their proof is to exploit previous results by \citeauthor{Seitz1974} \cite{Seitz1974} on the irreducible projective representations of simple groups of Lie type, later extended in \cite{Seitz1993} by \citeauthor{Seitz1993}. These results, alongside the well known projective representation theory of the alternating groups (see, for instance, \cite{Hoffman1992}) and of the sporadic groups \cite{ATLAS}, enable a description of all irreducible projective representations of degree $p$ of the simple groups (here we use the CFSG). 
	
	Discerning between $\rho(S)$ primitive or imprimitive, we list these groups in the following theorem.
	\begin{theorem}[\citeauthor{DixonZalesskii1998Primitive} \cite{DixonZalesskii1998Primitive}, \cite{DixonZalesskii2008Corrigendum}] \label{thm:socle_arb_p}
		Let ${\rho}(G)\leq \SL{p}{\C}$ be a finite primitive group with $p$ prime. Suppose $S \coloneqq \mathrm{soc}(G/Z(G))$ the socle of $G/Z(G)$ is a non-abelian simple group, then $G\leq Aut(S)$ and $S$ is one of the following:
		\begin{enumerate}
			\item $S$ is primitive and
			\begin{enumerate}
				\item $S\cong A_{p+1}$, with $p\geq 7$;
				\item $S\cong \mathrm{PSL}(2,q)$ where $p$ and $q$ satisfy one of the following conditions:
				\begin{enumerate}
					\item $p=q$ and $p\geq 11$,
					\item $p=\frac{q-1}{2}$, where $q$ a prime or $q=3^{\ell}$ with $\ell$ an odd prime,
					\item $p=\frac{q+1}{2}$, where $q=\ell^{2^k}\geq 5$ for some odd prime $\ell$ and $k\geq 0$,
					\item $p=2^\ell-1$ where $q=2^\ell$ for some odd prime $\ell$;
				\end{enumerate}
				\item $S\cong \mathrm{PSp}(2n,q)$ where $p,n$ and $q$ satisfy one of the following conditions:
				\begin{enumerate}
					\item $p=\frac{q^n+1}{2}$ where $n=2^{s}$ for $s\geq 2$, $q=\ell^{2^k}$ with $\ell$ an odd prime and $k\geq 0$,
					\item $p=\frac{3^{n}-1}{2}$ where $n$ is an odd prime and $q=3$;
				\end{enumerate}
				\item $S\cong \mathrm{PSU}(n,q)$ where $p=\frac{q^n+1}{q+1}$ and $n$ is an odd prime,
				\item $S$ satisfies one of the six exceptional cases:
				\begin{enumerate}
					\item $p=3$ with $S\cong \mathrm{PSL}(2,9)\cong A_6$,
					\item $p=7$ with $S\cong \mathrm{PSp}(6,2)$,
					\item $p=11$ with $S\cong M_{12}$,
					\item $p=23$ with $S\cong Co_2$, $Co_3$ or $M_{24}$.
				\end{enumerate} 
			\end{enumerate}
			\item $S$ is imprimitive, $G'$ is imprimitive and for some $n$ and $q$ we have that $G'\cong \mathrm{PSL}(n,q)$ and $p=\frac{q^n-1}{q-1}$. Furthermore, if $n=2$ then $q$ is even, otherwise $q$ is odd except for $(n,q) = (3,2)$.
		\end{enumerate}
		Conversely, whenever the parameters satisfy the suitable conditions, there is at least one  primitive group $H \leq \SL{p}{\C}$ (or conjugacy class) such that $\mathrm{soc}(H) \cong S$.
	\end{theorem}
	\begin{obs} 
		We may compare the groups described in this theorem with those in Theorem~\ref{thm:groups_small_p}. Note that the description given in Theorem~\ref{thm:groups_small_p} is more explicit. That is, instead of only providing the socles it provides the groups themselves. For instance, for $p=5$ there are two groups $G$ with socle $S=A_6$, these are $A_6$ and $S_6$. Recall that $G\leq Aut(S)$, in this case, $\mathrm{Aut}(A_6) = S_6$ and so these are all possible groups. Theorem~\ref{thm:groups_small_p} states that $S_6$ is in fact primitive. A priori, this last fact is not necessarily guaranteed by Theorem~\ref{thm:socle_arb_p} and Lemma~\ref{lemma:DZ_cases}.
	\end{obs}
	\begin{corrections}\label{correction:small_p_2}
		Following up on Correction~\ref{correction:small_p_1}, let us analyze Theorem~\ref{thm:socle_arb_p}(2) when $p=5,7,11$. If $p=5$, the only positive integer solutions to $5=\frac{q^n-1}{q-1}$ are $(n,q) = (5,1), (2,4)$. As $1$ is not a prime power, the only solution is $(2,4)$ and $S = \PSL{}{2,4} \cong A_5$. Note $A_5$ has a unique irreducible 5-dimensional projective representation (see \cite{ATLAS} or use \cite{GAP4}), and is given by a linear representation. This representation must be imprimitive by Theorem~\ref{thm:socle_arb_p}(2). Indeed, this has been confirmed both in Lemma~\ref{lemma:A_5_imprimitive} and by using \cite{Monomial2021}. Furthermore, the theorem states $\mathrm{Aut}(A_5) = S_5$ is primitive (there are no subgroups in between), which we already know from \ref{thm:groups_small_p} (3).
		
		When $p = 7$, the only pairs possible are $(n,q) = (2,6), (3,2)$. Since $6$ is not a prime power, we discard it. Furthermore, note $\PSL{}{2,7} \cong \PSL{}{3,2}$ and it has a unique irreducible $7$-dimensional projective representation, which is given by a linear representation. Again, by Theorem~\ref{thm:socle_arb_p}(2) we have that this group must me imprimitive. This has been confirmed using \cite{Monomial2021}.
		
		Finally, when $p=11$ no valid solutions $(n,q)$ can be found, so there is nothing to correct concerning this case. In further investigations, the situation of Theorem~\ref{thm:socle_arb_p}(2) should be carefully considered, as it seems to have been overlooked multiple times in the literature.
	\end{corrections}
	
	This concludes the known description of the groups corresponding to Lemma~\ref{lemma:DZ_cases}(ii). \\
	
	We move on to describe the case of Lemma~\ref{lemma:DZ_cases}(i), which also corresponds to Lemma~\ref{lemma:p_sylow_p}(ii) and (iii). An extensive study of this type of groups was done by \citeallauthors{KangZhang2009} in \cite{KangZhang2009}.
	
	\begin{lemma}[{\citeallauthors{KangZhang2009} \cite[Proposition 2.3]{KangZhang2009}}] \label{lemma:equiv_prim_monomial}
		Let $\rho(G) \leq \SL{p}{\C}$ be a finite primitive group. The following are equivalent:
		\begin{enumerate}[(i)]
			\item Let $S \coloneqq \mathrm{soc}(G/Z(G))$, then $S$ is an elementary abelian group of order $p^2$.
			\item There exists $H\trianglelefteq G$ such that $\rho(H)$ is imprimitive (monomial) and contains a non-scalar diagonal matrix.
		\end{enumerate}
	\end{lemma}
	
	Thus, the groups we wish to study are primitive groups with a normal imprimitive subgroup and with non-abelian socle. 
	\begin{obs}
		Recall that the primitive case of prime degree with a normal imprimitive subgroup and non-abelian socle was already covered in Theorem~\ref{thm:socle_arb_p}(2).	
	\end{obs}
	
	 For $p=3$, \citeauthor{Blichfeldt1917} in \cite[\S\S 78-79]{Blichfeldt1917} explicitly determined all primitive groups with a normal imprimitive subgroup. The technique employed by \citeauthor{Blichfeldt1917} in \cite{Blichfeldt1917} relied on a concept of ``invariant triangles'', however, this concept was only used in passing. 
	 
	 A different approach was later given by \citeauthor{Rigby1960} in \cite{Rigby1960}, to recover \citeauthor{Blichfeldt1917}'s results for $p=2,3$. Moreover, \citeauthor{Rigby1960} gave some interesting general results for arbitrary $p$, in particular, see \cite[Theorems 7, 8]{Rigby1960}. 
	 
	 More recently, \citeallauthors{KangZhang2009} generalized \citeauthor{Blichfeldt1917}'s concept of ``invariant triangles'' in \cite{KangZhang2009} to consider ``invariant polygons''. We shortly present their approach in the sequel.
	
	\begin{definition}
		Let $\rho(G) \leq \GL{n}{\C}$ be a finite imprimitive monomial group, define $V = \bigoplus_{j=0}^{n-1} \C e_j$ where $\{e_j\}_{j=0}^{n-1}$ is the canonical basis of $\C^n$. A $\rho(G)$-polygon $\Delta = \{v_0,\dots, v_{n-1}\}$ is a set of $n$ vectors of $V$ satisfying:
		\begin{enumerate}[(i)]
			\item $V = \sum_{j=0}^{n-1} \C v_j$,
			\item for any $g\in \rho(G)$ and for any $j=0,\dots, n-1$ we have $g(v_j)\in\C v_k$ for some $k=0,\dots,n-1$.
		\end{enumerate}
		Two $\rho(G)$-polygons $\Delta = \{v_j\}_{j=0}^{n-1},\, \Delta'= \{u_j\}_{j=0}^{n-1}$ are equal if there exist scalars $c_j\in\C^{\times}$ such that $v_j = c_j u_j$ for all $j=0,\dots,n-1$.
	\end{definition}
	\begin{obs}
		It can be seen that any finite monomial imprimitive group $\rho(G)$ has finitely many $\rho(G)$-polygons \cite[Lemma 3.1]{KangZhang2009}.
	\end{obs}
	For the convenience of the reader, we provide a summary of the relevant polygons and some elements of $\SL{p}{\C}$ used in the statements of the results.
	\begin{definition} \label{def:reps}
		Let $\{e_j\}_{j=0}^{p-1}$ be the canonical basis of $\C^p$. Let $\zeta = e^{2\pi i/p}$ be a primitive $p$-th root of unity. Define $\sigma, \tau,\lambda_d \in \SL{p}{\C}$ as follows
		\begin{equation*}
			\sigma: e_j\mapsto e_{j+1}, \quad \tau: e_j\mapsto \zeta^je_{j}, \quad \lambda_d: e_j \mapsto \epsilon e_{dj}
		\end{equation*}
		where $d\not\equiv 0 \mod p$, $0\leq j \leq p-1$, $\epsilon \in \C^{\times}$ is to be adjusted such that $\det(\lambda_d) = 1$. The $e_{j}$ with $j\geq p$ are $e_k$ with $0\leq k\leq p-1$ and $j\equiv k\mod p$. Furthermore, define the sets $\Delta_\infty = \{e_j\}_{j=0}^{p-1}$ and $\Delta_0,\dots, \Delta_{p-1}$ as follows:
		\begin{itemize}
			\item[] $\Delta_0 = \{u_j\}_{j=0}^{n-1}$ where $u_j = \sum_{k=0}^{p-1}\zeta^{jk}e_k$.
			\item[] $\Delta_i = \{\sigma^j(w_i)\}_{j=0}^{n-1}$ where $w_i = \sum_{k=0}^{p-1}\zeta^{i\binom{k}{2}}e_k$ for $1\leq i\leq p-1$.
		\end{itemize}
	\end{definition}
	
	With the introduced notation, we can proceed to state and describe the main results of \cite{KangZhang2009} which include \cite[Theorems 2.5, 2.6, 2.7]{KangZhang2009}. We start with a structural description of primitive groups with a normal monomial subgroup containing a non-scalar diagonal matrix.
	
	\begin{theorem}[\citeallauthors{KangZhang2009} \cite{KangZhang2009}] \label{thm:structure_primitive_normal_monomial}
		Let $p\neq 2$ be a prime. Let $\rho'(G) \leq \SL{p}{\C}$ be a finite primitive group with $H \trianglelefteq G$ such that $\rho'(H)$ is monomial containing a non-scalar diagonal matrix. Then there exists an equivalent representation $\rho:G \to \SL{p}{\C}$ such that 
		\begin{enumerate}[(i)]
			\item $\rho(D) = \langle\sigma, \tau\rangle \trianglelefteq \rho(G)$.
			\item $\Delta$ is a $\rho(D)$-polygon if and only if $\Delta\in\{\Delta_\infty,\Delta_0,\dots, \Delta_{p-1}\}$.
			\item $G$ acts on $\{\Delta_\infty,\Delta_0,\dots, \Delta_{p-1}\}$ where $g(\Delta_i)$  is defined by $\{\rho(g)(v_{ij})\}_{j=0}^{p-1}$, $\Delta_i = \{v_{ij}\}_{j=0}^{p-1}$ and $i\in \{\infty,0,\dots,p-1\}$ via $g(\Delta_i) = \Delta_{\tilde g(i)}, \, \tilde g \in S_p$.
			\item The group action of (iii) induces a non-trivial group morphism $\phi: G \to \PSL{}{2,p}$ with $\ker(\phi) = D$ or $\Ker(\phi) = \langle D, \lambda_{p-1}\rangle$ if $\lambda_{p-1}\notin \rho(G)$ or $\lambda_{p-1}\in G$ respectively. Furthermore, this morphism can be lifted to $\Phi: G \to \SL{}{2,p}$ with $\Ker(\Phi) = D$ and $\pi\circ\Phi = \phi$ where $\pi:\SL{}{2,p} \to \PSL{}{2,p}$ is the canonical projection.
			\item For any $g\in G$, if the action of $g$ is known, i.e $\phi(g)$ is known, then the representation $\rho(h)$ of some element $h\in G$ such that $g\in h\Ker(\phi)$ is known.
		\end{enumerate}
	\end{theorem}
	\begin{obs}
		The proof of Theorem~\ref{thm:structure_primitive_normal_monomial} can be found as follows: for  (i)-(iii) see \cite[Theorem 2.5]{KangZhang2009}(A)-(C). Part (iv) is a summary of \cite[Theorem 2.5]{KangZhang2009}(D) and \cite[Theorem 2.6]{KangZhang2009}. Part (v) corresponds to \cite[Theorem 2.5]{KangZhang2009}(E), whose statement is too long to fully reproduce here. However, it includes explicit representations of certain elements of $\rho(G)$ in terms of their action over $\{\Delta_\infty,\Delta_0,\dots, \Delta_{p-1}\}$. Thus, it is crucial to later determine explicit representations of these groups.
	\end{obs}
	\begin{obs}
		The group $D$ given by $\rho(D) = \langle\sigma, \tau\rangle$ in Theorem~\ref{thm:structure_primitive_normal_monomial}(i) is precisely the order $p^3$ group described in Remark~\ref{rmk:bis} and $D/Z(D)$ is the $p^2$ elementary abelian group in Lemma~\ref{lemma:equiv_prim_monomial}, corresponding to the socle of $G/Z(G)$.
	\end{obs}
	A rundown of the steps and ideas required for the proof of Theorem~\ref{thm:structure_primitive_normal_monomial} is given in more detail in \cite[p.1023]{KangZhang2009}. Nevertheless, we provide a couple of remarks that may be of interest to the reader.
	\begin{obs}
		In Theorem~\ref{thm:structure_primitive_normal_monomial}, the group $D$ corresponds to the extraspecial group of order $p^3$ described in Remark~\ref{rmk:case_i} and whose quotient by its center is the elementary abelian $p$-group described in Lemma~\ref{lemma:DZ_cases}(i). This alongside the action/morphism described in (iv) produce a semi-direct product as indicated by Lemma~\ref{lemma:DZ_cases}(i) when $p>2$.
	\end{obs}
	\begin{obs}
		The morphism $\phi$ arising from the action described in Theorem~\ref{thm:structure_primitive_normal_monomial}(iii) can be interpreted as follows. Consider the action by conjugation of $G$ over $D/Z(D)$. Since $D/Z(D) \cong \mathbb{F}_p\times\mathbb{F}_p$ can be interpreted as a $2$-dimensional $\mathbb{F}_p$-vector space, we get a representation of $G$ onto $\GL{}{2,p}$ and so projecting we get $\phi:G \to \PGL{}{2,p}$. Interpreting the set $\{\Delta_\infty,\Delta_0,\dots, \Delta_{p-1}\}$ as $\mathbb{P}^1(\mathbb{F}_p) = \{0,\dots,p-1,\infty\}$ it can be seen that this induces the same action on the $\rho(D)$-polygons as described in the theorem. 
	\end{obs}
	
	An application of this theorem leads to the following description of this type of groups.
	
	\begin{prop}[{\citeallauthors{KangZhang2009} \cite[Theorem 2.7]{KangZhang2009}}] \label{prop:pre_thm_primitive_mono}
		Let $p\neq 2$ be a prime. Define $f_1,f_2,f_3 \in \SL{p}{\C}$ as
		\begin{equation*}
			f_1: e_j \mapsto c_1\zeta^{\binom{j}{2}}e_j, \quad f_2:  e_j \mapsto c_2\sum_{k=0}^{p-1}\zeta^{jk}e_k, \quad f_3: e_j \mapsto c_3 e_{jm}
		\end{equation*}
		where $m\in\mathbb{F}^\times_p$ generates $\mathbb{F}^\times_p$ and $c_1,c_2,c_3\in \C^\times$ chosen appropriately for the $f_i$ to have determinant $1$. Let $\rho_0(G_0) = \langle\sigma,\tau,f_1,f_2,f_3\rangle \leq \SL{p}{\C}$, recall $\sigma, \tau$ from Definition~\ref{def:reps}. Then $\rho_0(G_0)$ is a primitive group of order $p^4(p^2-1)$ that contains a normal monomial subgroup with a non-scalar diagonal matrix.
	\end{prop}
	\begin{theorem}[{\citeallauthors{KangZhang2009} \cite[Theorem 2.7]{KangZhang2009}}] \label{thm:primitive_normal_monomial}
		Let $p\neq 2$ be a prime. Let $\rho'(G) \leq \SL{p}{\C}$ be a finite primitive group with $H \trianglelefteq G$ such that $\rho'(H)$ is monomial containing a non-scalar diagonal matrix. Then there exists an equivalent representation $\rho:G \to \SL{p}{\C}$ such that $\rho(D) \leq \rho(G) \leq \rho_0(G_0)$ and 
		\begin{enumerate}[(i)]
			\item If $p^4 \mid |G|$ then $G=G_0$. 
			\item If $p^4 \nmid |G|$ then $G$ is a semidirect product of $D$ and a subgroup of $\SL{}{2,p}$.
		\end{enumerate}
	\end{theorem}

	\begin{obs} \label{rmk:explicit_theorem}
		This theorem provides an explicit representation for case (i). Furthermore, in case (ii), note $|\SL{}{2,p}|=p(p^2-1)$ and $|D|=p^3$. Thus, if $p^4\nmid |G|$ this means that $G/D$ has order coprime to $p$ and is a subgroup of $\SL{}{2,p}$. These groups are known up to conjugation and can be found collected in \cite[Theorem 7.5]{KangZhang2009}. Furthermore, explicit generators for these subgroups in $\SL{}{2,p}$ are computed in \cite[Proposition 8.8, Theorems 8.9, 8.11, 8.13]{KangZhang2009}.
	\end{obs}
	
	Combining Theorem~\ref{thm:primitive_normal_monomial}, Remark~\ref{rmk:explicit_theorem} and Theorem~\ref{thm:structure_primitive_normal_monomial}(v), it is possible to obtain explicit generators for a representation of the groups in Theorem~\ref{thm:primitive_normal_monomial}(ii). \citeallauthors{KangZhang2009} computed these generators for $p=5,7$, making explicit the description of these groups given originally by \citeauthor{Brauer1967} in \cite{Brauer1967} for $p=5$ and by \citeauthor{Wales1969} in \cite{Wales1969},\cite{Wales1970}  for $p=7$. These generators can be found in \cite[Theorems A.3, A.6]{KangZhang2009}, moreover, these groups are gathered in the tables provided in Section~\ref{sec:tables} alongside additional information (see Observations below the table).
	
	We believe a computational implementation of the results in \cite{KangZhang2009} should be feasible, enabling the construction of representations for higher prime degrees, by following the ideas showed in the proofs of \cite[Theorems A.3, A.6]{KangZhang2009}.
	
	Finally, we conclude the post-CFSG approach subsection by recapping all results in the following classification statement.
	\begin{classification}[Primitive groups of prime degree] \label{class:primitive_any_p}
		Let $G$ be a finite group and $p\neq 2$ a prime. Let $\rho(G)\leq \SL{p}{\C}$ be a finite primitive group with $S \coloneqq \mathrm{soc}(G/Z(G))$. Then one of the following situations holds:
		\begin{enumerate}[(i)]
			\item $S$ is an elementary abelian $p$-group of order $p^2$ and $G/Z(G)$ is isomorphic to a subgroup $H$ of a semidirect product of $S$ and $\SL{}{2,p}$. A description of these groups is provided in Theorem~\ref{thm:primitive_normal_monomial} and Remark~\ref{rmk:explicit_theorem}. Structural results are covered in Theorem~\ref{thm:structure_primitive_normal_monomial}. A method for explicit computation of representations for these groups is exemplified for $p=5,7$ in \cite[Theorem A.3, A.6]{KangZhang2009}.
			\item $S$ is a non-abelian simple group and $G/Z(G)$ is isomorphic to a subgroup of $\mathrm{Aut}(S)$. Moreover, $\rho|_S$ is irreducible (either primitive or imprimitive). All possible socles $S$ in terms of the degree $p$ are listed in Theorem~\ref{thm:socle_arb_p}.
		\end{enumerate}
		This classification trivially translates to a classification for $\PGL{p}{\C}$ under the same hypotheses and we can replace $G/Z(G)$ by $G$ as $Z(G)=1$.
	\end{classification}
	
	\subsection{The composite dimensional case} \label{sec:composite}
	
	So far, most results and discussions have involved $n$, the degree of $\PGL{n}{\C}$, to be a prime number. In this subsection we explore results and classification statements when $n$ is composite. 
	
	Recall that the socle $S$ of a finite group is a direct product $S = G_1\times \cdots \times G_r$ of finite simple groups. So far, in the prime degree case, the only situation in which $r>1$ has been when the socle is abelian. Otherwise, $r=1$ is all we have seen in the non-abelian socle case. The following theorem shows the situation is different in the composite degree scenario.
	
	\begin{theorem} \label{thm:tensor_is_primitive}
		Let $\rho_i: G_i \to \GL{n_i}{\C}$, $i=1,2$ be (quasi-)primitive representations of $G_1$ and $G_2$. Then the Kronecker product of the representations $\rho_1 \otimes_\Kr \rho_2 : G_1\times G_2 \to \GL{n_1n_2}{\C}$ is a (quasi-)primitive representation of $G_1\times G_2$.
	\end{theorem}
	
		For a proof of Theorem~\ref{thm:tensor_is_primitive} see \cite[Theorem 1, Theorem 2]{Aschbacher2000}, where \citeauthor{Aschbacher2000} proves this result under more general hypotheses, see also Remark~\ref{rmk:aschbacher_general}.
		
	\begin{obs}
		The Kronecker product of representations is sometimes also called tensor product of representations. However, in some contexts this can be confused with the usual tensor product of representations, given $\rho_i: G\to \GL{n_i}{\C}$, $i=1,2$ the tensor product $\rho_1\otimes\rho_2 : G \to \GL{n_1n_2}{\C}$. To avoid confusion when $G_1 = G_2$, we use the $\otimes_{\Kr}$ notation.
	\end{obs}
	
	The first to explore the finite primitive groups of $\PGL{n}{\C}$ in a composite situation was \citeauthor{Blichfeldt1917}, who studied the case $n=4$ in \cite[VII]{Blichfeldt1917}. However, it was \citeauthor{Lin69} in his PhD thesis \cite{Lin69} who provided a general description of the groups present in the composite case.
	
	\begin{theorem}[{\citeauthor{Lin69} \cite[Chp III]{Lin69}}] \label{thm:composite_cases} Consider $K$ an algebraically closed field.
		Let $\rho: G\rightarrow GL_n(K)$ and the associated $\overline{\rho}: G/Z(G)\rightarrow \mathrm{PGL}_n(K)$ be faithful irreducible representations of degree $n$ of a finite group $G$ and assume that $(char(K),n)=1$
		. For a subgroup $\overline{H}$ of $G/Z(G)$, we denote by $H$ the inverse image of the projection map from $G$ to $G/Z(G)$.
		Then, one of the following items occur:
		\begin{enumerate}
			\item the representation is reducible or not quasi-primitive,
			\item the representation is a subgroup of a Kronecker product of two representations of smaller degree,
			\item $G/Z(G)$ has a normal abelian subgroup $\overline{P}$ of the form $C_p\times\ldots\times C_p$, the degree $n$ is $p^s$ with $p$ prime and $s$ natural, and $G/\overline{P}$ is isomorphic to a subgroup of the symplectic group $\mathrm{Sp}(2s,p)$. Moreover $\rho(P)$ modulo $Z(\mathrm{GL}_n(K))$ is up to conjugation by a fix matrix equal to $L\otimes_{\Kr}\ldots\otimes_{\Kr} L=:L_{\otimes^n_{\Kr}}$ linear group of order $p^{2s+1}$ and $L$ is a group in $GL_p(K)$ of order $p^3$ generated by the matrices $A=\mathrm{diag}(1,\xi,\ldots,\xi^{p-1})$ and $B$ a permutation matrix with $B^{-1}AB=\xi A$ where $\xi$ is a primitive $p$-th root of unity and as usual $\otimes_{\Kr}$ denotes the Kronecker product of matrices.  
			\item $G/Z(G)$ has a non-abelian normal subgroup $\overline{G}_1\times\cdots\times\overline{G}_m$ where $\overline{G}_i$ are simple groups which are conjugated and at least two primes divide $|G_i|$. Moreover the degree is $n=s^m$ for certain $s$, where $\rho_{|_{G_1}}(g)$ modulo conjugation corresponds to $I_m\otimes_{\Kr} \rho_1(g)$ with $\rho_1$ a irreducible representation of $G_1$ of degree $s$ and $I_m$ is the $m\times m$ identity matrix, and $\rho_{|\langle G_1,\ldots,G_m\rangle}$ modulo conjugation corresponds to $\otimes_{\Kr}^m \rho_1(G_1)$.
		\end{enumerate}
	\end{theorem}
	We have not found this result published outside of the PhD thesis \cite[Chp III]{Lin69}, so we provide a sketch of the proof for the reader's convenience.
	\begin{proof}[Proof (Sketch)] Consider $N$ a minimal normal subgroup of $G$ containing $Z(G)$, 
		by Clifford's Theorem (see, for example, \cite[\S1]{Clifford1937}), ${\rho}_{|N}$ is irreducible or a direct sum of equivalent irreducible representations of $N$. If all constituents are not identical, then the representation is not quasi-primitive or is reducible \cite[Lemma 51.2]{Curtis1962}. Assume that all such conjugate representations are identical.
		When ${\rho}_{|N}$ in not irreducible, by \cite[\S3]{Clifford1937} (see also \cite[\S2]{Huppert1957}) we have $$\rho_{|N}(N)=I_r\otimes_{\Kr}\mathfrak{N}=\left(\begin{array}{cccc}
			\mathfrak{N}&0&\ldots&0\\
			0&\mathfrak{N}&\ldots&0\\
			\vdots&\vdots&&\vdots\\
			0&0&\ldots&\mathfrak{N}\\
		\end{array}\right)$$
		with $\mathfrak{N}$ a finite subgroup of $\GL{k}{K}$. By \cite[Satz 3]{Huppert1957} (see also \cite[Theorem (51.7)]{Curtis1962}), there exist representations $\rho_1, \rho_2$ of $G$ of degrees $k$ and $n/k$ respectively such that ${\rho}(g)={\rho}_1(g)\otimes_{\Kr} {\rho}_2(g)$ for all $g\in G$.
		
		Now consider ${\rho}_{|N}$ is irreducible. Write $\overline{N}=\overline{G}_1\times\ldots\times \overline{G}_m$ where the $\overline{G}_i$ are simple groups which are conjugated. Assume $|G_i|=p^k$ with $p$ a prime and $k\geq 1$. Since the degree of a representation divides the index of any abelian subgroup (of $N$ in our situation), we obtain that $n\mid p^s$ and the socle in such case is an elementary abelian group. The details and description presented in the statement follow from \cite[Lemma1,Lemma2,Chapter III]{Lin69} which we do not reproduce here.
		
		Now assume at least two different primes $p$ and $q$ divide $|\overline{G}_i|$.	\citeauthor{Lin69} observed that the centralizer of $G_1$ contains $N$ since the $G_i$ themselves commute element-wise. That is, if $g_1\in G_1$ and $g_j\in G_j$ then $g_1^{-1}g_jg_1 g_j^{-1}=(g_1,g_j)\in Z(G)$, and if $g_1$ has order $p$ and $g_j$ has order $q$ then $(g_1,g_j)$ has order $p$ and $q$, thus trivial. 
		
		Next, the argument in \cite{Lin69} is as follows.
		Once such commutativity is fixed and $N$ is inside the subgroup generated by $G_1$ and $C_G(G_1)$, the centralizer of $G_1$ by elements of $G$, then $\rho_{|G_1}=I_r\otimes_{\Kr} {\rho_1}_{|G_1}$ where $\rho_1$ is an irreducible representation of $G_1$ and $r$ a natural. Then, by \cite[Satz 3]{Huppert1957} $\rho_{|\langle G_2,\ldots,G_m\rangle}=\rho'_{|\langle G_2,\ldots,G_m\rangle}\otimes_{\Kr} I_s$ with $s$ the degree of $\rho_1$ and $\rho'$ a linear group representation for $\langle G_2,\ldots,G_m\rangle$. Now, because $G_2$ is conjugate to $G_1$, $\rho_{|G_2}=I_t\otimes_{\Kr} {\rho_1}_{|G_1}\otimes_{\Kr} I_s$ for some natural $t$, and iterating such procedure one can conclude up to conjugation that $$\rho_{|\langle G_1,\ldots,G_m\rangle}=I_u\otimes_{\Kr}{\rho_1}_{|G_1}\otimes_{\Kr}\ldots\otimes_{\Kr} {\rho_1}_{|G_1}$$
		but because $\rho$ restricted to $N$ is irreducible, we have that $u=1$, and the order of the representation is $s^m$.
	\end{proof}
	
	\begin{obs}
		If $n$ is a prime number, then a quasi-primitive group must fall either in case 3 or 4. Indeed, case 3 recovers Lemma~\ref{lemma:DZ_cases}(i), since $s=1$ and $\mathrm{Sp}(2,p) \cong \SL{}{2,p}$. Case 4 recovers Lemma~\ref{lemma:DZ_cases}(ii) since $n=s^m$ implies $m=1$.
	\end{obs}
	In order to provide some examples of each situation, we introduce the classifications of $\SL{n}{\C}$ for $n=4,6$ proved by \citeauthor{Blichfeldt1917} in \cite[VII]{Blichfeldt1917} and \citeauthor{Lin69} in \cite{Lin71} respectively.
	\begin{theorem}[{\citeauthor{Blichfeldt1917} \cite[VII]{Blichfeldt1917}}] \label{thm:n_4}
		Let $\rho(G) \leq \SL{4}{\C}$ be a finite primitive group, then one of the following is satisfied
		\begin{enumerate}[(i)]
			\item $G$ is quasisimple and $G/Z(G)\cong A_5, A_6, A_7, \PSL{}{2,7}, \PSU{}{4,2}$.
			\item $G/Z(G) \cong S_5, S_6$.
			\item $G/Z(G) \cong A\times B$ where $A, B \cong S_4, A_4, A_5$ for any possible combination or an extension  $[G/Z(G):A\times B] = 2$ if $A=B$ or $[G/Z(G):A_4\times A_4] = 4$.
			\item $G/Z(G)$ has a normal abelian subgroup $P = C_2^4$ and $(G/Z(G))/P \cong C_5, D_5, \mathrm{Sz}(2), A_5, S_5,$ $A_6, S_6$.
		\end{enumerate} 
	\end{theorem}
	
	\begin{obs}
		All the described groups may be found with greater detail in Table~\ref{tab:pgl_4}. Furthermore, \citeauthor{Blichfeldt1917} provided the explicit representations for each group, these can be found in \cite[VII]{Blichfeldt1917} using the identifiers given in Table~\ref{tab:pgl_4}.
	\end{obs}
	
	\begin{example}
		The groups in Theorem~\ref{thm:n_4}(i),(ii) correspond to case Theorem~\ref{thm:composite_cases}(4) with $s=4$ and $m=1$. Moreover, $A_5\times A_5$ and its extension also fall in this case.
		
		For $n=4$ and case Theorem~\ref{thm:composite_cases}(3), necessarily $s=2, p=2$ and so $\mathrm{Sp}(4,2) \cong S_6$. The groups corresponding to this situation are the products and extensions involving $A_4$ and $S_4$, since they have a normal abelian subgroup $C_2^4$ and their quotients are subgroups of $S_6$. Furthermore, all groups in Theorem~\ref{thm:n_4}(iv) are also of this type. 
		
		Case Theorem~\ref{thm:composite_cases}(2) includes all $A \times B$ for $A, B \cong S_4, A_4, A_5$, which partially overlaps with the previous situation.
	\end{example}
	
	\begin{theorem}[\citeauthor{Lin71} \cite{Lin71}] \label{thm:n_6}
		Let $\rho(G) \leq \SL{6}{\C}$ be a finite quasi-primitive group, then one of the following is satisfied
		\begin{enumerate}[(i)]
			\item $G/Z(G) \cong A \times B$ where $A \cong A_5$, $A_4$, or $S_4$, and $B \cong \mathrm{PSL}(2, 7)$, $A_5$, $A_6$, or $H_i$ for $i = 1, 2$, or $3$, where $H_3$ is the Hessian group isomorphic to an extension of $C_3 \times C_3$ by $\mathrm{SL}(2, 3)$, and $H_2$ and $H_1$ have indices $3$ and $6$ in $H_3$.
			
			\item $G/Z(G)$ is isomorphic to a subgroup of index $2$ in $S_4 \times H_i$, $i = 1$ or $2$, or to a subgroup of index $3$, $12$, or $24$ in $A_4 \times H_3$.
			
			\item $G$ is quasisimple and $G/Z(G)\cong A_5, A_6, A_7, \PSL{}{2,7}, \PSL{}{2,11}, \PSL{}{2,13},$ $\PSU{}{4,2}$, $\PSU{}{3,3}$, $\PSU{}{4,3}, \PSL{}{3,4}$ or $J_2$ the Hall-Janko group.
			
			\item $G/Z(G) \cong S_5, S_7$ or $[G/Z(G):H]=2$ with $H \cong A_6, \PSL{}{3,4}, \PSU{}{3,3}, \PSU{}{4,2}$ is a split extension by $C_2$.
		\end{enumerate}
	\end{theorem}
	\begin{obs}
		All the described groups in items (iii) and (iv) may be found with greater detail in Table~\ref{tab:pgl_6}. Note that this description is not fully explicit. In item (ii), there might be multiple non-isomorphic groups of a certain index within the group. Thus, it would remain to make explicit how many there are and their description. Still, this is just routine computations.
	\end{obs}
	\begin{obs}
		In the case $n=6$ note that no groups of type Theorem~\ref{thm:composite_cases}(3) can appear since $6$ is not a prime power.
	\end{obs}
	We proceed to introduce the known results for $n=8,9,10$. In the case $n=8$ the classification is essentially known up to precise determination of groups coming from tensor products. The classification is due to \citeauthor{Huffman1976} in \cite{Huffman1976} and \citeauthor{Feit1976}, however, Feit's results relied largely on computational work and seem to be unpublished, only an outline of the proof is known to be found in \cite[Theorem A]{Feit1976}.
	
	For $n=9$, the case where $7\mid |G|$ was partially described in \citeauthor{Doro1976}'s PhD thesis \cite{Doro1976} (supervised by Feit) and later refined by Feit. The groups are fully covered in \cite[Theorem B]{Feit1976}. The case when $7\nmid |G|$ is known to be groups with $|G| = 2^a3^b5^c$ for some $a,b,c \geq 0$, but lack any further description as far as we know.
	
	When $n=10$, not much seems to be known. \citeauthor{Feit1976} in \cite{Feit1976} gives a brief statement on this case, and \citeauthor{Doro1976}'s PhD thesis contains some groups with $7\mid |G/Z(G)|$ which potentially appear when $n=10$. However, \citeauthor{Doro1976} does not distinguish between $n=9$ and $n=10$, so it is to be seen which groups correspond to which case.
	
	\begin{theorem}[{Huffman, Feit and Wales \cite{Huffman1976}, \cite[Theorem A]{Feit1976}}] \label{thm:n_8}
		Let $\rho(G) \leq \SL{8}{\C}$ be a finite quasi-primitive group, then one of the following is satisfied
		\begin{enumerate}[(i)]
			\item $G$ is quasisimple and $G/Z(G) \cong A_6, A_8, A_9$, $\PSL{}{2,7}$, $\PSL{}{2,17}$, $\PSL{}{2,8}$, $\mathrm{Sp}(2,6)$, $\mathrm{PS}\Omega^+(8,2)$.
			\item $G/Z(G) \cong S_8, S_9$, $\PGL{2}{7}$ or an extension $[G/Z(G):H]=2$ with $H \cong A_6, \mathrm{PS}\Omega^+(8,2)$ or an extension $[G/Z(G):H]=3$ with $H \cong \PSL{}{2,8}$.
			\item $G/Z(G) \cong A\times B$ with $A, B$ groups with a projective faithful quasi-primitive representations of degree $2$ and $4$.
			\item $O_2(G/Z(G)) \neq 1$ and $G/O_2(G) \cong A \leq \mathrm{Sp}(6,2)$ where $O_2(G)$ is the maximal $2$-group of $G$.
			\item $G/Z(G)$ is isomorphic to an extension of $A_5^3$ by $C_3$ or $S_3$.
		\end{enumerate}
	\end{theorem}
	\begin{obs}
		In the statement of Theorem~\ref{thm:n_8}, $\mathrm{PS}\Omega^+(8,2)$ refers to the simple group $D_4(2)$ in Lie type group notation. The original statement introduces the group $\mathrm{W}(\mathrm{E}_8)$, the Weil group of type $\mathrm{E}_8$. In particular, the groups appearing are $\mathrm{W}(\mathrm{E}_8)'$ and $\mathrm{W}(\mathrm{E}_8)$ itself. Quotient the center, the first is the simple group and the second the extension. Sometimes $\mathrm{O}^+(8,2)$ is also used, but this notation may refer to a different group, so we avoid it. 
	\end{obs}	
	\begin{theorem}[{\citeauthor{Doro1976}, \citeauthor{Feit1976} \cite{Doro1976}, \cite[Theorem B]{Feit1976}}] \label{thm:n_9}
		Let $\rho(G) \leq \SL{9}{\C}$ be a finite quasi-primitive group, then one of the following is satisfied
		\begin{enumerate}
			\item $|G| = 2^a3^b5^c$ for some $a,b,c\geq 0$.
			\item $G$ is quasisimple and $G/Z(G) \cong A_{10}, \PSL{}{2,7}$, $\PSL{}{2,19}$.
			\item $G/Z(G) \cong S_{10}$ or $G/Z(G) \cong \PSL{}{2,7} \times A$ for some $A$ group with a projective faithful quasi-primitive representations of degree $3$ or $G/Z(G)$ is an extension of degree $2$ of $\PSL{}{2,7} \times \PSL{}{2,7}$.
		\end{enumerate}
	\end{theorem}
	
	We collect the presented results in the following classification statement.
	\begin{classification}[Quasi-primitive groups of small composite degree] \label{class:primitive_small_composite}
		Let $4\leq n\leq 10$ be a composite number. Consider $\rho(G) \leq \SL{n}{\C}$ a finite irreducible group. 
		\begin{enumerate}[(a)]
			\item If $n=4$, \citeauthor{Blichfeldt1917} classified the primitive groups $\rho(G)$ in \cite[VII]{Blichfeldt1917}. These groups are covered in Theorem~\ref{thm:n_4} and can be consulted in Table~\ref{tab:pgl_4}. The explicit representations of all groups are available in \cite[VII]{Blichfeldt1917}.
			\item If $n=6$, \citeauthor{Lin69} classified all quasi-primitive groups $\rho(G)$ in \cite{Lin69}. These are covered in Theorem~\ref{thm:n_6} and can be found in Table~\ref{tab:pgl_6}. An explicit description of all groups originating from tensor products is missing.
			\item If $n=8$, \citeauthor{Feit1971}, \citeauthor{Huffman1976} classified all quasi-primitive groups $\rho(G)$ in \cite{Huffman1976} and \cite{Feit1976}. However, \citeauthor{Feit1976}'s complete proof seems to be unpublished. The groups are covered in Theorem~\ref{thm:n_8}.
			\item If $n=9$, \citeauthor{Feit1976} and \citeauthor{Doro1976} partially classified the quasi-primitive groups $\rho(G)$ in \cite{Feit1976} and \cite{Doro1976}. The only explicitly known groups are those such that $7$ divides $|G|$. The groups are covered in Theorem~\ref{thm:n_9}.
			\item If $n=10$, partial results exist on the quasi-primitive groups $\rho(G)$ in \cite{Feit1976} and \cite{Doro1976}. However, a complete classification is not known to the author.
		\end{enumerate}
	\end{classification}
	
  	\newpage
  	\section{Quasisimple groups and computational implementations} \label{sec:quasisimple_computational}
  	
  	\subsection{Some results concerning quasisimple groups} \label{subsec:quasisimple}
  	
  	The focus of this survey has been on the description of all finite subgroups of $\GL{n}{\C}$ and $\PGL{n}{\C}$. Nevertheless, there is a class of groups in which this theory is greatly more developed than the general case: quasisimple groups. Recall a group $G$ is quasisimple if $G'=G$ and it is a central extension of a simple group. Furthermore, any covering of a simple group is quasisimple, so these give rise to all projective irreducible representation of the simple groups. Finally, also note that for simple groups, quasi-primitive and irreducible are equivalent notions.
  	
  	In this section, we describe some of the more relevant results regarding linear and projective irreducible representations of simple and quasisimple groups. Many of the presented results are either given in the form of formulas and tables. As such, the access to this results is very tedious when done by hand. This is the main motivation for Section~\ref{subsec:computational}, where we present a computer program based on Python whose main objective is to improve accessibility to many of the results in the sequel. See  Section~\ref{subsec:computational} and \cite{FISGO} for more information.
  	
  	We start by presenting a famous result of the late 90's due to \citeauthor{Tiep1996} \cite{Tiep1996}. This result gives a description of all irreducible representations of relatively small degree of the quasisimple groups.
  	
  	\begin{theorem}[\citeauthor{Tiep1996} \cite{Tiep1996}] \label{teoTZ96}
  		Let $G$ be a quasisimple irreducible complex linear group of degree $d$. Assume that
  		$d\leq2r$ for some prime divisor $r$ of $|G|$, and let $L=G/Z(G)$. Then one of the following holds.
  		\begin{enumerate}[itemsep=0pt, topsep=3pt,noitemsep]
  			\item $L=A_n$, $\max\{9, r\}\leq n \leq 2r+1$, $d=n-1$.
  			\item $L = \PSL{}{2,q}$, $q\neq 5,7,9$, and one of the following holds:
  			\begin{enumerate}[itemsep=0pt, topsep=0pt,noitemsep]
  				\item $q =2^a$, $a\geq 3$, $r=2^a\pm 1$ is a Fermat or a Mersenne prime, $d\in\{r,r \mp 1, r \mp 2\}$;
  				\item  $q =r\geq11$, $d\in\{r, (r \pm 1)/2, r \pm 1\}$;
  				\item  either $q\geq 11$ is a prime or $q=3^n$, with $n$ an odd prime. Furthermore, $r = (q-1)/2$ and $d\in\{ r , r+1 , 2r\}$.
  				\item either $q\geq 13$ is a prime or $q=5^n$, $n$ an odd prime. Furthermore, $r=(q-1)/4$ and $d=2r$;
  				\item $q\geq 13$, $r=(q+1)/2$, and $d\in\{r-1,r, 2r-2 , 2r-1 , 2r \}$ ;
  				\item either $q\geq 11$ is a prime or $q=3^n$, $n$ an odd prime. Furthermore, $r =(q+1)/4$ and $d\in\{2r-1, 2r\}$.
  			\end{enumerate}
  			\item $L=\PSL{}{n,q}$, $n\geq 3$, and one of the following holds:
  			\begin{enumerate}[itemsep=0pt, topsep=0pt,noitemsep]
  				\item $q = 2$, $n\geq 5$, either $r=2^{n-1}-1$ or $r =2^n-1$, and $d =2n-2$ (1 rep);
  				\item $q\geq 3$, $n$ an odd prime, $r = (q^n-1)/(q-1)$. Furthermore, $d =r-1$ (1 rep) or $d = r$
  				($q-2$ reps);
  			\end{enumerate}
  			\item $L= \PSU{}{n,q}$, $n\geq 3$, and one of the following holds:
  			\begin{enumerate}[itemsep=0pt, topsep=0pt,noitemsep]
  				\item  $q = 2$, $n-1\geq 5$ is an odd prime, $r = (2^{n-1} + 1)/3$. Furthermore, $d = 2r - 1$ (2 reps) or
  				$d = 2r$ (1 rep);
  				\item $n$ is an odd prime, $r = (q^n+ 1)/(q + 1)$. Furthermore, $d = r - 1$ (1 rep) or $d = r$ ($q$
  				reps).
  			\end{enumerate}
  			\item $L=\mathrm{PSp}(2n,q)$, $n\geq 2$, and one of the following holds:
  			\begin{enumerate}[itemsep=0pt, topsep=0pt,noitemsep]
  				\item $q = 3$, $n$ an odd prime, $r = (3^n-1)/2$. Furthermore, $d = r$ (2 reps) or $d = r + 1$ (2
  				reps);
  				\item $q = 3$, $n$ an odd prime, $r = (3^n + 1)/4$. Furthermore, $d = 2r - 2$ (2 reps) or $d = 2r$ (2
  				reps);
  				\item  $q = 5$, $n an odd prime$, $r = (5^n-1)/4$. Furthermore, $d = 2r$ (2 reps);
  				\item  $n = 2^m$, $r = (q^n+1)/2$. Furthermore, $d = r-1$ (2 reps) or $d = r$ (2 reps).
  			\end{enumerate}
  			\item  Exceptions for alternating and finite classical groups:
  			\begin{enumerate}[itemsep=0pt, topsep=0pt,noitemsep]
  				\item $L = A_5 = \PSL{}{2,4} = \PSL{}{2,5}$, $(r, d) = (2, 2)$, $(2, 3)$, $(2, 4)$, $(3, 2)$, $(3, 3)$, $(3, 4)$, $(3, 5)$, $(3, 6)$, $(5, 2)$,
  				$(5, 3)$, $(5, 4)$, $(5, 5)$, $(5, 6)$;
  				\item $L = A_6 = \PSL{}{2,9} = \Sp{}{4,2}'$, $(r,d) = (2,3)$, $(2,4)$, $(3,3)$, $(3,4)$, $(3,5)$, $(3,6)$, $(5,3)$, $(5,4)$, $(5,5)$,
  				$(5, 6)$, $(5, 8)$, $(5, 9)$, $(5,10)$; 
  				\item $L = \PSL{}{3,2} = \PSL{}{2,7}$, $(r,d) = (2,3)$, $(2,4)$, $(3,3)$, $(3,4)$, $(3,6)$, $(7,3)$, $(7,4)$, $(7,6)$, $(7,7)$, $(7,8)$;
  				\item $L = \SL{}{3,3}$, $r = 13$, $d = 12$ (1 rep), $d = 13$ (1 rep), $d=16$ (4 reps)
  				or $d = 26$ (3 reps);
  				\item $L = A_7$, $(r, d) = (2, 4)$, $(3, 4)$, $(3, 6)$, $(5, 4)$, $(5, 6)$, $(5, 10)$, $(7, 4)$, $(7, 6)$, $(7, 10)$, $(7, 14)$;
  				\item $L = \PSL{}{3,4}$, $(r, d) = (3, 6)$, $(5, 6)$, $(5, 8)$, $(5, 10)$, $(7, 6)$, $(7, 8)$, $(7, 10)$;
  				\item $L = A_8 = \SL{}{4,2}$, $(r, d)= (5, 7)$, $(5, 8)$, $(7, 7)$, $(7, 8)$, $(7, 14)$;
  				\item $L=A_9$, $r = 7$, $d = 8$ (2 reps);
  				\item $L = A_{11}$, $r = 11$, $d = 16$ (2 reps);
  				\item $L = \PSL{}{4,3}$, $r = 13$, $d = 26$ (2 reps);
  				\item $L = \SU{}{3,3}$, $(r, d) = (3, 6)$, $(7, 6)$, $(7, 7)$, $(7, 14)$;
  				\item $L=\SU{}{4,2} = \PSp{}{4,3}$, $(r,d)= (2,4)$, $(3,4)$, $(3,5)$, $(3,6)$, $(5,4)$, $(5,5)$, $(5,6)$, $(5,10)$;
  				\item $L = \PSU{}{4,3}$, $r = 3, 5, 7$, $d = 6$ (4 reps);
  				\item $L = \SU{}{5,2}$, $r = 5$, $d = 10$ (1 reps);
  				\item $L = \Sp{}{6,2}$, $r = 5$, $7$, $d = 7$ (1 reps) and $d = 8$ (1 reps);
  				\item $L = \Sp{}{4,4}$, $r = 17$, $d = 18$ (1 reps) or $d = 34$ (2 reps);
  				\item $L=\Omega_8^+(2)$, $r = 5$, $7$, $d = 8$ (1 reps);
  				\item $L=\Omega_8^{-}(2)$, $r= 17$, $d = 34$ (1 reps).
  			\end{enumerate}
  			\item $L$ is an exceptional group of Lie type:
  			\begin{enumerate}[itemsep=0pt, topsep=0pt,noitemsep]
  				\item $L={}^2{\!B}_2(8)$, $r = 7$, $13$, $d = 14$ (2 reps);
  				\item $L={}^3{\!D}_4(2)$, $r = 13$, $d = 26$ (1 rep.);
  				\item $L=G_2(3)$, $r = 7$, $13$, $d = 14$ (1 rep);
  				\item $L=G_2(4)$, $r = 7$, $13$, $d = 12$ (1 rep);
  				\item ${}^2{\!F}_4(2)'$ , $r = 13$, $d = 26$ (2 reps).
  			\end{enumerate}
  			\item $L$ is a sporadic finite simple group:
  			\begin{enumerate}[itemsep=0pt, topsep=0pt,noitemsep]
  				\item $L = M_{11}$ either $r = 5$, $11$ and $d : 10$ (3 reps), or 
  				$r=11$ and $d=11$ (1 rep),
  				or $r=11$ and $d=16$ (2 reps);
  				\item $L=M_{12}$, either $r= 5$, $11$ and $d = 10$ (2 reps), or 
  				$r = 11$ and $d = 11$ (2 reps),
  				or $r=11$ and $d = 12$ (1 rep), or $r = 11$ and $d=16$ (2 reps);
  				\item  $L=M_{22}$, either $r = 5$, $7$, $11$ and $d = 10$ (2 reps), or $r=11$ and $d = 21$ (3 reps);
  				\item $L =J_2$ either $r= 3$, $5$, $7$ and $d=6$ (2 reps), or $r = 7$ and $d=14$ (3 reps);
  				\item $L = M_{23}$, either $r = 11$, $23$ and $d = 22$ (1 rep), or 
  				$r = 23$ and $d = 45$ (2 reps);
  				\item $L = HS$, $r = 11$ and $d = 22$ (1 rep);
  				\item $L = J_3$, $r = 17$, $19$ and $d = 18$ (4 reps);
  				\item $L = M_{24}$, $r = 23$, $d = 23$ (1 rep) or $d = 45$ (2 reps);
  				\item  $L = McL$, $r=11$ and $d = 22$ (1 rep);
  				\item $L = Ru$, $r = 29$ and $d = 28$ (2 reps);
  				\item $L = Suz$, $r = 7$, $11$, $13$ and $d = 12$ (2 reps);
  				\item $L = Co_3$, $r = 23$ and $d = 23$ (1 rep);
  				\item $L = Co_2$, $r = 23$ and $d = 23$ (1 rep);
  				\item $L = Co_1$, $r = 13$, $23$ and $d = 24$ (1 rep).
  			\end{enumerate}
  		\end{enumerate}
  		Conversely, if the triple  $(L,r,d)$ satisfies any of these conditions, then some covering group $G$ of $L$ has an irreducible complex representation of dimension $d$.
  	\end{theorem}
  	\begin{obs}
  		This theorem includes the representations of the socles described in Theorem~\ref{thm:socle_arb_p}, so it can be considered as a generalization. However, here we do not distinguish between primitive and imprimitive representations, nor are limited to considering a prime degree.
  	\end{obs}
  	
  	It is interesting to consider separately the ingredients used to prove Theorem~\ref{teoTZ96}. The first main ingredient includes knowledge on the first few degrees of the irreducible representations of the classical groups SL, SU, Sp. The second ingredient is a complete collection of the minimal degree of the projective irreducible representations of the simple groups of Lie type, excluding the trivial representation.
  	
  	\begin{classification}[Minimal degree of an irreducible projective simple group of Lie type in zero characteristic] \label{class:small_pirreps}
  		Let $G(q)$ denote a finite simple group of Lie type defined over a field of order $q = p^k$ for some prime $p$ and integer $k\geq 1$. Let $d(G(q))$ denote the smallest integer $n > 1$ such that $G(q)$ has a projective irreducible representation of degree $n$ over $\C$. Then $d(G(q))$ is known and was collected by \citeauthor{Tiep2000} in \cite[Table 1, Table 2]{Tiep2000}. The results are originally due to \citeauthor{Tiep1996} in \cite{Tiep1996} for the classical groups, and \citeauthor{Lubeck2001} in \cite{Lubeck2001} for the exceptional groups. 
  	\end{classification}
  	
  	The previous classification only concerns the zero characteristic case. Nevertheless, there is also much information known concerning the cross characteristic case. We refer the interested reader to \cite[\S 4, \S 5]{Tiep2000}, \cite{Seitz1974} and \cite{Seitz1993} among others, we apologize for the many omissions.
  	
  	These results culminated in a classification, due to \citeauthor{HissMalle2001} \cite{HissMalle2001, HissMalle2002}, of all (absolutely) irreducible linear representations of the quasisimple groups in cross-characteristic of degree less than $250$, with some omissions complemented by \citeauthor{Lubeck2001} in \cite{Lubeck2001Complement}. In the zero characteristic case, this is equivalent to the description of all projective irreducible representations of the simple groups of degree less than $250$.
  	
  	\begin{classification}(Absolutely irreducible linear representations of degree less than $250$ of the quasisimple groups in cross-characteristic) \label{class:HissMalle}
  		\begin{itemize}
  			\item[(a)] \citeauthor{HissMalle2001} provide in \cite[Table 2]{HissMalle2001} and \cite[Table 2]{HissMalle2002} all the groups in the classification statement, alongside additional information on each representation, except for the representations of the Lie groups in their defining characteristic.
  			\item[(b)] \citeauthor{Lubeck2001} in \cite{Lubeck2001Complement} describes all absolutely irreducible representations of the quasisimple finite Lie groups in their defining characteristic up to degree $l^3/8$ where $l$ is the rank of the Lie group. These include all those of degree at most 250.
  		\end{itemize}
  	\end{classification}
  	
  	Finally, we highlight some computational results due to \citeauthor{Lubeck2001} which can be found in his website \cite{LuebeckWebsite}. These results concern all degrees of irreducible complex representations, together with their multiplicities, of some groups of Lie type of rank at most $8$. These include, among others, the $\SL{l+1}{q}$, $\SU{l+1}{q}$, $\Sp{2l}{q}$ and $\mathrm{Spin}_{2l+1}(q), \mathrm{Spin}_{2l}(q)$ groups for the classical groups; the exceptional groups ${}^3D_4(q), {}^2E_6(q), E_6(q), E_7(q), E_8(q)$ and the Suzuki and Ree groups ${}^2B_2(q^2), G_2(q), {}^2G_2(q^2), F_4(q), {}^2F_4(q^2)$; where $l$ denotes the rank. The results are given for general $q$ a prime power and in the format of GAP formulas. All data can be accessed in \cite{LuebeckWebsite}.

  	\subsection{Computational implementations} \label{subsec:computational}
  	
  	
  	The statement of Theorem~\ref{teoTZ96}, \citeauthor{HissMalle2001}'s tables in Classification~\ref{class:HissMalle} containing more than 1000 entries, henceforth the ``Hiss-Malle table''; and the large amount of representation data provided by \citeauthor{Lubeck2001} in \cite{LuebeckWebsite} are all a testament to how tedious it can be to systematically access all this important representation information.
  	
  	While it is true that GAP \cite{GAP4} has a large library of character tables, and is very useful to work with the characters and groups themselves, it is hard to answer a question like ``\textit{Which simple groups have projective irreducible representations of degree 150?}''. This is a question which can be answered with ease using Classification~\ref{class:HissMalle}. However, there is no guarantee that all groups present in the Hiss-Malle table are actually stored in GAP, for not all of their characters may have ever been computed, only the small ones might actually be known.
  	
  	An alternative software providing some solutions to this problem is CHEVIE \cite{Chevie}, a program split between GAP \cite{GAP4} and Maple \cite{Maple}, focused on working with Lie type groups. The Maple part of CHEVIE implements some of the generic character tables (character tables generated from one or more parameters) concerning the groups of \cite{LuebeckWebsite}. While it does not specifically implement all the data provided in \cite{LuebeckWebsite}, it solves the part where, if the character table of a group is not known, then it can be generated.
  	
  	Nevertheless, this is still an overkill approach to the type of question which we now propose.
  	
  	\begin{question}
  		For a fixed integer $n \geq 2$, what simple groups have degree $n$ projective irreducible representations over $\C$?
  	\end{question}
  	
  	To help answer this specific question, which is a subset of the main question concerning this survey, the author has set out to collect all generic information on representation degrees of the simple groups and collect it in a single package for the mathematical community to access. Although still in active development by the author, we would like to present some of the thought process behind the program and describe its main functionalities. 
  	
  	Firstly, the program is written in pure Python, an will be distributed as a Python package in the future. This has the advantage of being independent of a specific software other than the base programming language itself. Furthermore, it can easily be used alongside SageMath \cite{SageMath}, which is written in C and Python. Since SageMath \cite{SageMath} provides an interface to other systems such as GAP \cite{GAP4} and Maple \cite{Maple}, it might even be feasible to interface with CHEVIE \cite{Chevie}. Furthermore, much of the available data is provided as precomputed information, and Python is excellent to deal with such data.
  	
  	We called this program FiSGO: Finite simple groups by order \cite{FISGO}. Its name is due to it being part of an earlier project, which was expanded into this survey. All information on the program, functionality, status of the development and documentation can be found in
  	
  	\begin{center}
  		\url{https://github.com/GeraGC/FiSGO}
  	\end{center}
  	
  	We proceed to list some of the main functionalities and capabilities of FiSGO \cite{FISGO}.
  	
  	\begin{enumerate}[(1)]
  		\item Given a list of prime factors $N = 2^a3^b5^c\cdots$ for $a,b,c, \dots \geq 0$ integers, identify all simple groups whose order divides $N$. It is also possible to specify an absolute bound, such as $M = 10!$, so that only groups of order less than $M$ and dividing $N$ are returned.
  		\item Functions implementing \citeauthor{Blichfeldt1917} and \citeauthor{Brauer1967}'s theorems \ref{thm:bound_primes}, \ref{thm:bound_primes_ged_dim}, \ref{thm:bound_primes_order} and \ref{thm:prime_bounds_general} bounding the prime powers of quasi-primitive groups. An implementation of \citeauthor{Collins2007}' absolute bound in Theorem~\ref{thm:collins_bounds} is also provided.
  		\item Simple group objects are available for all families of simple groups. These provide access to basic information such as the order of the group, its Schur multiplier, recommended notations in LaTeX, its GAP/ATLAS notation, among other capabilities. 
  		\item A system of ``\textit{simple group codes}'' has been implemented. Each code is a string of words and numbers identifying a specific simple group. For example, \verb|CA-1-13| refers to the Chevalley group of type A with parameters $n=1, q = 13$, also known as $\PSL{2}{13}$. This provides a memory efficient and unified way to identify, store and access information on a specific group simply by knowing its code.
  		\item An (almost complete) interface to the extended Hiss-Malle table, containing all information stored in \cite[Table 2]{HissMalle2002} alongside the omissions specified in \cite[Table 2]{HissMalle2001}. The only data missing is the ``field'' column in the case of the omissions specified in \cite[Table 2]{HissMalle2001}. \citeauthor{Lubeck2001}'s data \cite{Lubeck2001Complement} on representations in the defining characteristics is to be eventually included completing all information on cross-characteristic absolutely irreducible representations of degree less than 250 of the quasi-simple groups.
  		\item All information on ordinary character degrees and multiplicities of the sporadic groups and their coverings.
  		\item Access to the minimal degree of a complex irreducible projective representation for any simple group as stated in Classification~\ref{class:small_pirreps}. Eventually, the bounds on the cross-characteristic case may be also implemented for easier access.
  		\item An interface to all data on provided by \citeauthor{Lubeck2001} in \cite{LuebeckWebsite} concerning all degrees of irreducible complex representations, together with their multiplicities, of the (non-exceptional) covering groups of Lie type of rank at most $8$. Data on groups whose Schur cover is exceptional, meaning it is larger than expected, is also manually taken from the GAP database.
  		\item An implementation of \citeauthor{Tiep1996}'s Theorem~\ref{teoTZ96} concerning the relatively small complex irreducible representations of the quasisimple groups.
  		\item A function to search for all complex projective irreducible representation data of the simple groups stored in the program (as files or formulas) for a fixed degree or a range of degrees. This includes a list of all groups with a representation in the requested range of degrees (alongside its corresponding degrees). The function also provides two additional lists, one containing all groups whose displayed information is complete, meaning there are no other representations in the chosen range apart from the ones provided. The second list contains all groups whose information is partial, meaning other representations in the chosen range of degrees could exist.
  	\end{enumerate}
  	
  	Originally, FiSGO started as the an implementation of points (1), (2) and (10) to search for primitive simple groups, and so the origin of the name. From the stated points, currently points (8) and (9) are in active developement, the rest are already available. Moreover, the Hiss-Malle tables are available in various formats (csv, ASCII txt, Markdown,...) in the GitHub repository, so anyone can easily access and treat the information itself.
  	
  	All documentation is available in its separate website, accessible from GitHub \cite{FISGO}. Documentation is automatically produced as the code is written and implemented, so it is readily available at any moment. A more detailed list of progress and an introductory tutorial may be found on GitHub \cite{FISGO}. The origin of all used data is also stated and credited in the documentation website, the functions docstrings and the GitHub readmes.\\
  	
  	As an ending note, we add the fact that, once points (8) and (9) are successfully implemented, it should be relatively easy to produce a table containing for each degree $2 \leq n \leq N$, all the simple groups that have a complex projective irreducible representation for $N$ at least $1000$. This list would be omitting the alternating groups, whose representation results are not yet planned to be implemented, but could be obtained with other software.

	\newpage
	\section{Tables of primitive subgroups of $\PGL{n}{\C}$ for $2\leq n\leq 11$} \label{sec:tables}
	
	\subsection{How to read the tables} \label{subsec:how_to_tables}
	\begin{itemize}
		\item In column ``$G$'' a standard name for the group is provided. If the given group does not have a standard/known name, then a generic identifier is provided.
		\item In column ``$|G|$'' a description of the order of the group is given.
		\item In column ``Origin'' denotes the origin of the representations. ``C'' denotes that the projective representations originate from a covering, ``L'' denotes them originating from linear representations of the group itself.
		\item The \cite{GAP4} column provides an indentifier (if possible) for the group in GAP's database of groups. There are mainly two databases used: the Small Group library (S) and the Perfect Groups library (P). In the observations column the letters ``S'' and ``P'' denote which GAP database has been used for each identifier.
		\item The StructureDescription column contains the output (or expected output) of the function \verb|StructureDescription| from \cite{GAP4}. This function is intended to provide a \textbf{qualitative} description of the group structure. It is not an isomorphism invariant, and isomorphic groups may even produce different results. For detailed information on the symbols used and the output produced consult the \cite{GAP4} documentation.
		\item In column ``FI'' the number of non-isomorphic faithful irreducible representations of the group (of the appropriate dimension) is provided, if known. These correspond to the linear representations inducing the projective representations. It is possible that some of the indicated representations are projectively equivalent, meaning they induce the same projective representation. The data has been computed using \cite{GAP4} character tables.
		\item In the \emph{Obs} column, additional information is provided for a given group. The letters ``S'' and ``P'' denote which GAP database has been used for each identifier in the GAP column.
		\item For $n=4$, the table contains an additional column labeled \cite{Blichfeldt1917}. This column denotes the name of each group as given in \citeauthor{Blichfeldt1917}'s book \cite[Chapter VII]{Blichfeldt1917}.
	\end{itemize}
	\textbf{Additional comments}
	\begin{itemize}
		\item Some of the groups organized by blocks in Table~\ref{tab:pgl_4} for $\PGL{4}{\C}$ are related among each others. The reader interested in subgroup lattices involving these groups may consult them in \cite[Annex A]{Cheltsov2019} or \cite[\S 4]{Avila2024}. Furthermore, both \cite[Annex A]{Cheltsov2019} and \cite[\S 5]{Avila2024} provide the generators of the representation of the described groups in a more collected format than \cite[VII]{Blichfeldt1917}.
		\item Many of the naming conventions used in this tables have been taken from the database GroupNames \cite{GroupNames} by \citeauthor{GroupNames}. It has also been used alongside GAP \cite{GAP4} to verify properties of the groups described in the tables.
	\end{itemize}
	\newpage
	
	\subsection{Primitive subgroup tables of $\PGL{n}{\C}$ for $2\leq n \leq 7$} \label{subsec:tables}
	\setlength\intextsep{0mm}
	\renewcommand{\tablename}{\textbf{Table}}
	\renewcommand{\thetable}{\textbf{\arabic{table}}}
	\begin{table}[h!]
		\centering
		\caption{\bfseries Primitive subgroups of $\PGL{2}{\C}$}
		\label{tab:pgl_2}
		\renewcommand{\arraystretch}{1.25}
		\begin{tabular}{|l|l|c|l|c|c|c|} \hline
			$G$ & $|G|$ & Origin & GAP\cite{GAP4} & StructureDescription & FI & Obs \\ \hline\hline
			$S_4$ & $24$ & C & [24, 12] & \verb|S4| & 2 & S  \\ \hline
			$\PSL{2}{3} = A_4$ & $12$ & C & [12, 3] & \verb|PSL(2,3)| & 3 & S \\ \hline
			$A_5$ & $60$ & C & [60, 5] & \verb|PSL(2,5)| & 2 & S \\ \hline
		\end{tabular}
	\end{table}
	
	\begin{table}[h!]
		\centering
		\caption{\bfseries Primitive subgroups of $\PGL{3}{\C}$}
		\label{tab:pgl_3}
		\renewcommand{\arraystretch}{1.25}
		\begin{tabular}{|l|l|c|l|c|c|c|} \hline
			$G$ & $|G|$ & Origin & GAP\cite{GAP4} & StructureDescription & FI & Obs \\ \hline\hline
			$A_5$ & $60$ & L & [60,5] & \verb|A5| & 2 & S \\ \hline
			$\PSL{2}{7}$ & $168$ & L & [168,42] & \verb|PSL(3,2)| & 2 & S \\ \hline
			$A_6$ & $360$ & C & [360, 118] & \verb|A6| & 4 & S \\ \hline
			$H_1$ & $36$ & C & [36, 9] & \verb|(C3xC3):C4| & 8 & S, \ref{item:P3_2} \\ \hline 
			$H_2$ & $72$ & C & [72, 41] & \verb|(C3xC3):Q8| & 8 & S, \ref{item:P3_2} \\ \hline 
			$H_3$ & $216$ & C & [216, 153] & \verb|((C3xC3):Q8):C3| & 6 & S, \ref{item:P3_1}, \ref{item:P3_2} \\ \hline 
		\end{tabular}
	\end{table}
	\vspace{5mm}
	\textbf{Observations}
	\begin{enumerate}
		\item Group $H_3$ is known as the \emph{Hessian} group, also known as the affine special linear group $ASL_2(3) = SA(2,3)$. The given structure description refers to $H_3 = \mathrm{PSU}_3(2) \rtimes C_3$ and $\mathrm{PSU}_3(2) = C_3^2\rtimes Q_8$. An alternative structure description would be $H_3 = C_3^2\rtimes \SL{2}{3}$, which is coherent with observation \ref{item:P3_2}. \label{item:P3_1}
		\item The groups $H_i$ for $i=1,2,3$ are primitive groups containing a normal monomial subgroup $D\cong C_3\times C_3$. The groups $H_i$ are semidirect products of $D$ and a subgroup $K_i$ of $\SL{2}{3}$. In particular, $K_1 = C_4$, $ K_2 = Q_8$ the quaternion group and $K_3 = \SL{2}{3}$ itself. Generators for these groups can be found in \cite[\S 79]{Blichfeldt1917}.
		\label{item:P3_2}
	\end{enumerate} 
	\newpage
	

	\begin{table}[h!]
		\centering
		\caption{\bfseries Primitive subgroups of $\PGL{4}{\C}$}
		\label{tab:pgl_4}
		\renewcommand{\arraystretch}{1.25}
		\begin{tabular}{|l|c|l|c|l|c|c|} \hline
			$G$ & \cite{Blichfeldt1917} & $|G|$ & Origin & GAP\cite{GAP4} & StructureDescription & Obs \\ \hline\hline
			$A_5$ & (A),(B) & 60 & C,L & [60,5] & \verb|A5| &  S, \ref{item:P4_1} \\ \hline
			$A_6$ & (C) & 360 & C & [360,118] & \verb|A6| &  S, \ref{item:P4_1} \\ \hline
			$A_7$ & (D) & 2520 & C & [2520,1] & \verb|A7| &  P, \ref{item:P4_1} \\ \hline
			$\PSL{2}{7}$ & (E) & 168 & C & [168,42] & \verb|PSL(3,2)|  & S, \ref{item:P4_1} \\ \hline
			$\mathrm{PSU}_4(2)$ & (F) & 25920 & C & [25920,1] & \verb|O(5,3)|  & P, \ref{item:P4_1} \\ \hline
			$S_5$ & (G),(H) & 120 & C,C & [120,34] & \verb|S5| &  S, \ref{item:P4_1} \\ \hline
			$S_6$ & (K) & 720 & C & [720,763] & \verb|S6| &  S, \ref{item:P4_1} \\ \hline \hline
			
			$A_4^2$ & 1º & 144 & C & [144,184] & \verb|A4xA4| &  S, \ref{item:P4_2} \\ \hline
			$\mathrm{PSO}^+_4(3)$ & 2º & 288 & C & [288,1026] & \verb|(A4xA4):C2| &  S, \ref{item:P4_2} \\ \hline
			$A_4\times S_4$ & 3º & 288 & C & [288,1024] & \verb|A4xS4| &  S, \ref{item:P4_2} \\ \hline
			$A_5\times A_4$ & 4º & 720 & C & [720,768] & \verb|A5xA4| &  S, \ref{item:P4_2} \\ \hline
			$S_4\times S_4$ & 5º & 576 & C & [576,8653] & \verb|S4xS4| &  S, \ref{item:P4_2} \\ \hline
			$S_4\times A_5$ & 6º & 1440 & C & [1440,5848] & \verb|S4xA5| &  S, \ref{item:P4_2} \\ \hline
			$A_5^2$ & 7º & 3600 & C & [3600,1] & \verb|A5xA5| &  P, \ref{item:P4_2} \\ \hline \hline
			
			$\mathrm{PSO}^+_4(3)\rtimes C_2$ & 8º & 576 & C & [576,8654] & \verb|((A4xA4):C2):C2| &  S, \ref{item:P4_3} \\ \hline
			$A_4^2\rtimes C_4$ & 9º & 576 & C & [576,8652] & \verb|(A4xA4):C4| &  S, \ref{item:P4_3} \\ \hline
			$A_4\wr C_2$ & 10º & 288 & C & [288,1025] & \verb|(A4xA4):2| &   S, \ref{item:P4_3} \\ \hline
			$A_5\wr C_2$ & 11º & 7200 & C & ID1 & \verb|(A5xA5):2| &  \ref{item:P4_3} \\ \hline
			$S_4\wr C_2$ & 12º & 1152 & C & [1152,157849] & \verb|(S4xS4):2| &  S, \ref{item:P4_3} \\ \hline \hline
			
			$C_2^4\rtimes C_5$ & 13º & 80 & C & [80,49] & \verb|(C2xC2xC2xC2):C5| &  S, \ref{item:P4_4} \\ \hline
			$C_2^4\rtimes D_5$ & 14º & 160 & C & [160,234] & \verb|((C2xC2xC2xC2):C5):C2| &   S, \ref{item:P4_4} \\ \hline
			$C_2^4\rtimes \mathrm{Sz}(2)$ & 15º & 320 & C & [320,1635] & \verb|((C2xC2xC2xC2):C5):C4| &  S, \ref{item:P4_4} \\ \hline
			$C_2^4\rtimes A_5$ & 16º & 960 & C & [960,11358] & \verb|(C2xC2xC2xC2):A5| &   S, \ref{item:P4_4} \\ \hline
			$C_2^4\rtimes A_5$ & 17º & 960 & C & [960,11357] & \verb|(C2xC2xC2xC2):A5| &   S, \ref{item:P4_4} \\ \hline
			$C_2^4\rtimes S_5$ & 18º & 1920 & C & [1920,240996] & \verb|(C2xC2xC2xC2):S5| &  S, \ref{item:P4_4} \\ \hline
			$C_2^4\rtimes S_5$ & 19º & 1920 & C & [1920,240993] & \verb|(C2xC2xC2xC2):S5| &  S, \ref{item:P4_4} \\ \hline
			$C_2^4\rtimes A_6$ & 20º & 5760 & C & [5760,1] & \verb|(C2xC2xC2xC2):A6| &   P, \ref{item:P4_4} \\ \hline
			$C_2^4.S_6$ & 21º & 11520 & C &  & \verb|((C2xC2xC2xC2):A6):C2| &  \ref{item:P4_4} \\ \hline
		\end{tabular}
	\end{table}
	\begin{minted}[frame=lines,framesep=2mm,bgcolor=pastel3,linenos]{gap-console}
gap> ID1 := WreathProduct(AlternatingGroup(5),Group((1,2)));;
	\end{minted}
	\textbf{Observations}
	\begin{enumerate}
		\item Primitive simple group or group with a primitive simple normal subgroup.
		\label{item:P4_1}
		\item Primitive group with a normal reducible subgroup.
		\label{item:P4_2}
		\item Primitive group with one of the groups in \ref{item:P4_2} as a normal subgroup.
		\label{item:P4_3}
		\item Primitive group with a normal imprimitive subgroup.
		\label{item:P4_4}
	\end{enumerate}
	
	\newpage
	\begin{table}[h!]
		\centering
		\caption{\bfseries Primitive subgroups of $\PGL{5}{\C}$}
		\label{tab:pgl_5}
		\renewcommand{\arraystretch}{1.25}
		\begin{tabular}{|l|l|c|l|c|c|c|} \hline
			$G$ & $|G|$ & Origin & GAP\cite{GAP4} & StructureDescription & FI & Obs \\ \hline\hline
			$A_6$ & $360$ & L & [360,118] & \verb|A6| & 2 & S \\ \hline
						$S_5$ & $120$ & L & [120,34] & \verb|S5| & 2 & S \\ \hline
			$S_6$ & $720$ & L & [720,763] & \verb|S6| & 4 & S \\ \hline
			$\PSL{2}{11}$ & $660$ & L & [660,13] & \verb|PSL(2,11)| & 2 & S \\ \hline
			$\mathrm{PSU}_4(2)$ & $2^63^45$ & L & [25920,1] & \verb|O(5,3)| & 2 & P \\ \hline
			$G_1$ & $75$ & C & [75,2] & \verb|(C5xC5):C3| & 12 & S, \ref{item:P5_1} \\ \hline
			$G_2$ & $150$ & C & [150,6] & \verb|(C5xC5):C6| & 24 & S, \ref{item:P5_1} \\ \hline
			$G_3$ & $200$ & C & [200,44] & \verb|(C5xC5):Q8| & 16 & S, \ref{item:P5_1} \\ \hline
			$G_4$ & $300$ & C & [300,23] & \verb|(C5xC5):(C3:C4)| & 16 & S, \ref{item:P5_1} \\ \hline
			$G_5$ & $600$ & C & [600,150] & \verb|(C5xC5):SL(2,3)| &  & S, \ref{item:P5_1} \\ \hline
			$G_6$ & $3000$ & C & [3000,1] & \verb|(C5xC5):SL(2,5)| &  & P, \ref{item:P5_1} \\ \hline
		\end{tabular}
	\end{table}
	\textbf{Observations}
	\begin{enumerate}
		\item The groups $G_i$ for $i=1,\dots, 6$ are primitive groups containing a normal monomial subgroup $D \cong C_5\times C_5$. The groups $G_i$ are semidirect products of $D$ and a subgroup $K_i$ of $\SL{2}{5}$. In particular, $K_1 = C_3, \, K_2 = C_6, \, K_3=Q_8$ the quaternion group, $K_4 = \mathrm{Dic}_3 = C_3\rtimes C_4$ the dicyclic group of order $12$, $K_5 = \SL{2}{3}$ and $K_6 = \SL{2}{5}$. Generators for these groups can be found in \cite[Theorem A.3]{KangZhang2009} with the same numbering. 
		\label{item:P5_1}
		\item The original classification due to \citeauthor{Brauer1967} in \cite{Brauer1967} included the group $A_5$. However, it has been determined that this group is imprimitive, so it has been removed. See Corrections~\ref{correction:small_p_1}, \ref{correction:small_p_2}.
	\end{enumerate}
	
	\newpage
	\begin{table}[h!]
		\centering
		\caption{\bfseries Quasi-primitive subgroups of $\PGL{6}{\C}$}
		\label{tab:pgl_6}
		\renewcommand{\arraystretch}{1.25}
		\begin{tabular}{|l|l|c|l|c|c|c|} \hline
			$G$ & $|G|$ & Origin & GAP\cite{GAP4} & StructureDescription & FI & Obs \\ \hline\hline
			TP & - & - & - & - & - & \ref{item:P6_2} \\ \hline \hline
			$A_5$ & $60$ & C & [60,5] & \verb|A5| & 1 & S \\ \hline
			$S_5$ & $120$ & C & [120,34] & \verb|S5| & 2 & S \\ \hline
			$A_6$ & $360$ & C,C & [360,118] & \verb|A6| & 2,4 & S \\ \hline
			$A_7$ & $2520$ & L,C,C & [2520,1] & \verb|A7| & 1,2,4 & P \\ \hline
			$S_7$ & $5040$ & L & ID1 & \verb|S7| & 2 & \\ \hline
			$\PSL{2}{7} $ & $168$ & L,C & [168,42] & \verb|PSL(3,2)| & 1,2 & S \\ \hline
			$\PGL{2}{7} $ & $336$ & L,C & [336,208] & \verb|PSL(3,2):C2| & 3,8 & S \\ \hline
			$\PSL{2}{11}$ & $660$ & C & [660,13] & \verb|PSL(2,11)| & 2 & S \\ \hline
			$\PSL{2}{13}$ & $1092$ & C & [1092,25] & \verb|PSL(2,13)| & 2 & S \\ \hline
			$\PSU{4}{2}$ & $2^63^45^1 $ & L & [25920,1] & \verb|O(5,3)| & 1 & P \\ \hline
			$\PSU{3}{3}$ & $6048$ & L & [6048,1] & \verb|PSU(3,3)| & 1 & P \\ \hline
			$\PSU{4}{3}$ & $2^73^65^17^1$ & C & ID2 & \verb|PSU(4,3)| & 2 & \\ \hline
			$J_2$ & $2^73^35^27^1 $ & C & ID3 & \verb|HJ| & 2 & \\ \hline
			$\PSL{3}{4}$ & $2^63^25^17^1 $ & C & [20160,5] & \verb|PSL(3,4)| & 2 & P \\ \hline \hline
			$A_6.C_2$ & $720 $ & C & [720,765] & \verb|A6.C2| & 6 & S \\ \hline
			$\PSU{4}{2}\rtimes C_2$ & $2^73^45^1$ & L & ID5 & \verb|O(5,3):C2| & 2 & \\ \hline
			$\PSU{3}{3} \rtimes C_2$ & $12096$ & L & ID4 & \verb|PSU(3,3):C2| & 2 & \\ \hline
			$\PSL{3}{4} \rtimes C_2$ & $2^83^65^17^1$ & C,C & ID6 & \verb|PSL(3,4):C2| & 4,4 & \ref{item:P6_1} \\ \hline
		\end{tabular}
	\end{table}
	\begin{minted}[frame=lines,framesep=2mm,bgcolor=pastel3,linenos,fontsize=\footnotesize]{gap-console}
gap> ID1 := SymmetricGroup(7);;
gap> ID2 := ProjectiveSpecialUnitaryGroup(4,3);;
gap> ID3 := SimpleGroup("J2");;
gap> ID4 := Group([(3,4,7,15,18,26,21,14)(5,10,12,22,24,16,23,8)(6,13,11,9,19,20,17,25)(27,28),
(1,2,3,5,11,20,23,26)(4,8,18,21,25,28,6,14)(7,16,24,12,10,15,9,19)(22,27),(3,6)(4,9)(5,12)(7,17)
(10,16)(11,21)(13,15)(14,20)(18,19)(23,24)(25,26)(27,28)]);
gap> ID5 := Group([(2,4,7)(5,8,13)(6,12,9)(10,16,23)(11,19,15)(14,21,30)(17,22,31)(18,26,25)
(20,28,36)(24,32,33)(27,35,37)(29,34,40),(1,2,5,10,17,12)(3,6,11,18,21,7)(4,8,14)(9,15,22)
(13,20,29,36,30,38)(16,24,33,25,34,40)(19,27,32,35,31,39)(23,26),(1,3)(2,6)(4,9)(5,11)(7,12)
(8,15)(10,18)(13,19)(14,22)(16,25)(17,21)(20,27)(23,26)(24,34)(28,37)(29,32)(30,31)(33,40)
(35,36)(38,39)]);
gap> ID6 := Group([(1,42,12,24,22,32,23,40)(2,21,35,7,34,9,33,31)(3,17,37,26,41,18,30,11)
(4,16,36,14,38,19,6,10)(5,27,39,15)(8,28,20,29)(13,25),(1,23)(2,33)(3,30)(4,39)(5,38)(8,20)
(9,19)(10,21)(13,15)(14,18)(16,17)(24,29)(26,31)(28,40)(34,41)(35,37),(3,4)(5,6)(7,10)(8,9)
(11,15,14,12)(13,27,19,26)(16,23,22,31)(17,21,20,18)(24,39,36,40)(25,32,37,28)(29,42,33,41)
(30,38,34,35)]);
	\end{minted}
	
	\textbf{Observations}
	\begin{enumerate}
		\item These refer to groups arising from tensor products of representations of lower dimensions, a qualitative description can be found in \ref{thm:n_6}. \label{item:P6_2}
		\item Has 2 non-isomorphic covers with $|Z|=6$. \label{item:P6_1}
	\end{enumerate}
	\newpage

	\begin{table}[h!]
		\centering
		\caption{\bfseries Primitive subgroups of $\PGL{7}{\C}$}
		\label{tab:pgl_7}
		\renewcommand{\arraystretch}{1.25}
		\begin{tabular}{|l|l|c|l|c|c|c|} \hline
			$G$ & $|G|$ & Origin & GAP\cite{GAP4} & StructureDescription & FI & Obs \\ \hline\hline
			$A_8$ & $2^6 3^2 5^1 7^1$ & L & [20160,4] & \verb|A8| & 1 & P \\ \hline
			$S_8$ & $2^7 3^2 5^1 7^1$ & L & ID1 & \verb|S8| & 2 &  \\ \hline
			$\PSL{2}{13}$ & $1092$ & L & [1092,25] & \verb|PSL(2,13)| & 2 & S \\ \hline
			$\PGL{2}{7}$ & $336$ & L & [336,208] & \verb|PSL(3,2):C2| & 2 & S  \\ \hline
			$\mathrm{PSp}_6(2)$ & $2^9 3^4 5^1 7^1$ & L & [1451520,1] & \verb|O(7,2)| & 1 & P \\ \hline
			$\mathrm{PSU}_3(3)$ & $2^5  3^3  7^1$ & L & [6048,1] & \verb|PSU(3,3)| & 3 & P \\ \hline
			$\mathrm{G}_2(2)$ & $2^6  3^3  7^1$ & L & ID2 & \verb|U3(3).2| & 2 &  \\ \hline
			$\PSL{2}{8}$ & $504$ & L & [504,156]  & \verb|PSL(2,8)| & 4 & S \\ \hline
			$\mathrm{R}(3)$ & $1512$ & L & [1512,779] & \verb|PSL(2,8):C3| & 3 & S \\ \hline
			$G_1$ & $196$ & C & [196,8] & \verb|(C7xC7):C4| & 24 & S, \ref{item:P7_1} \\ \hline
			$G_2$ & $392$ & C & [392,36] & \verb|(C7xC7):C8| &  & \ref{item:P7_1}  \\ \hline
			$G_3, G_4$ & $392$ & C & [392,38] & \verb|(C7xC7):Q8| &  & \ref{item:P7_1} \\ \hline
			$G_5$ & $588$ & C & [588,33] & \verb|(C7xC7):(C3:C4)| &  & \ref{item:P7_1} \\ \hline
			$G_6$ & $784$ & C & [784,162] & \verb|(C7xC7):Q16| &  & \ref{item:P7_1} \\ \hline
			$G_7, G_8$ & $1176$ & C & [1176,215] & \verb|(C7xC7):SL(2,3)| &  & \ref{item:P7_1} \\ \hline
			$G_9, G_{10}$ & $2352$ & C &  & \verb|(C7xC7):(C2.S4)| &  & \ref{item:P7_1}, \ref{item:P7_2} \\ \hline
			$G_{11}$ & $16464$ & C & [16464,1] & \verb|(C7xC7):SL(2,7)| &  & P, \ref{item:P7_1}, \ref{item:P7_2} \\ \hline
		\end{tabular}
	\end{table}
	\begin{minted}[frame=lines,framesep=2mm,bgcolor=pastel3,linenos]{gap-console}
gap> ID1 := SymmetricGroup(8);;
gap> ID2 := Group(Group(AtlasGenerators("U3(3).2",2).generators));;
	\end{minted}
	
	\textbf{Observations}
	\begin{enumerate}
		\item The groups $G_i$ for $i=1,\dots, 11$ are primitive groups containing a normal monomial subgroup $D \cong C_7\times C_7$. The groups $G_i$ are semidirect products of $D$ and a subgroup $K_i$ of $\SL{2}{7}$. In particular, $K_1 = C_4, \, K_2 = C_8, \, K_3=K_4 = Q_8$ two non conjugate copies of the quaternion group within $\SL{2}{7}$, $K_5 = \mathrm{Dic}_3 = C_3\rtimes C_4$ the dicyclic group of order $12$, $K_6 = Q_{16}$ the generalized quaternion group, $K_7=K_8=\SL{2}{3}$ corresponding to two non-conjugate copies of $\SL{2}{3}$ within $\SL{2}{7}$, $K_9=K_{10} = S_4^+ = \mathrm{CSU}_2(3)$ two non-conjugate copies of the conformal special unitary group of dimension 2 over $\mathbb{F}_3$ within $\SL{2}{7}$, and finally $K_{11} = \SL{2}{7}$ itself. Generators for these groups can be found in \cite[Theorem A.6]{KangZhang2009} with the same numbering. 
		\label{item:P7_1}
		\item Groups calculated using Magma. In particular, $G_{11}$ has been identified as the unique perfect group of order 16464 by checking perfection.
		\label{item:P7_2}
		\item The original classification due to \citeauthor{Wales1970} in \cite{Wales1970} included the group $\PSL{}{2,7}$. However, it has been determined that this group is imprimitive, so it has been removed. See Corrections~\ref{correction:small_p_1}, \ref{correction:small_p_2}.
	\end{enumerate}
	\newpage
	\subsection{Missing classifications for $\PGL{n}{\C}$ for $2\leq n\leq 11$}
	In this final subsection we describe what is missing to accomplish a fully explicit classification of the (quasi)-primitive groups of $\PGL{n}{\C}$ for $2\leq n\leq 11$.
	
	We consider a classification \emph{complete} if a table, as the ones previously presented, can be given detailing all primitive groups for such $n$. This is a very restrictive criterion, as classifications of $n=6$ and $n=11$ are generally considered complete in the literature. \\
	
	\begin{table}[h!]
		\centering 
		\caption{\bfseries Status of the classification for small dimension}
		\label{tab:missing}
		\vspace{1mm}
		\renewcommand{\arraystretch}{1.25}
		\begin{tabularx}\hsize{|c|c|X|}
			 \hline
			$n$ & Complete & Comments / \tabitem Missing  \\ \hline\hline
			$2$ & Yes &  - \\ \hline
			$3$ & Yes &  - \\ \hline	
			$4$ & Yes &  \citeauthor{Blichfeldt1917} classified the primitive groups, not quasi-primitive \\ \hline
			$5$ & Yes &  - \\ \hline
			$6$ & No &  \makecell[{{>{\parindent0em}p{\hsize}}}]{
			Parts missing from Theorem~\ref{thm:n_6}: \\
			\tabitem An explicit description of all quasi-primitive groups arising from subgroups of tensor products and their explicit representations.} \\ \hline
			$7$ & Yes &  - \\ \hline
			$8$ & No & \makecell[{{>{\parindent0em}p{\hsize}}}]{
				Parts missing from Theorem~\ref{thm:n_8}: \\
				\tabitem Explicit description of all quasi-primitive groups arising from subgroups of tensor products and their explicit representations. \\
				\tabitem Explicit description of the group extensions of $A_6$, $\mathrm{PS}\Omega^+(8,2)$, $\PSL{}{2,8}$ and $A_5^3$.} \\ \hline
			$9$ & No &  \makecell[{{>{\parindent0em}p{\hsize}}}]{
				Parts missing from Theorem~\ref{thm:n_9}: \\
				\tabitem Explicit description of the groups $\PSL{}{2,7} \times A$ with $A$ quasi-primitive of degree $3$ \\
				\tabitem Explicit description of the groups $|G| = 2^a3^b5^c$. \\
				\tabitem Explicit description of the group extension of $\PSL{}{2,7}^2$
				}\\ \hline
			$10$ & No &  \tabitem Essentially unclassified \\ \hline
			$11$ & No &  \makecell[{{>{\parindent0em}p{\hsize}}}]{
				\tabitem From $G'$ given in Theorem~\ref{thm:groups_small_p}, description of the full groups $G$.
				\\
				\tabitem A list of all groups in $I_{11}$ (Lemma~\ref{lemma:equiv_prim_monomial} with $p=11$), up to conjugates in $G_0$ (Proposotion~\ref{prop:pre_thm_primitive_mono}) could be obtained alongside a representation using the results of \cite{KangZhang2009}, see Classification~\ref{class:primitive_any_p}(i).} \\
				 \hline
		\end{tabularx}
	\end{table}

	\clearpage
	\printbibliography
	
	\noindent Gerard Gonzalo Calbetó, \\
	
	\noindent Departament Matemàtiques, Edif. C, \\
	
	\noindent Universitat Autònoma de Barcelona, \\
	
	\noindent 08193 Bellaterra, Catalonia, Spain\\
	
	\noindent ggonzalo.math@gmail.com // gerard.gonzalo@uab.cat
	
\end{document}